\documentclass[11pt,a4paper]{article}
\usepackage[top=1in, bottom=1.25in, left=1in, right=1in]{geometry}
\usepackage[utf8]{inputenc}
\usepackage[english]{babel}
\usepackage{graphicx,epstopdf}
\usepackage{caption, subcaption,float}
\usepackage{color}
\usepackage{multirow}
\usepackage{tikz}
\usepackage{amsmath}
\usepackage{amsthm}
\usepackage{tabls}
\usepackage{csquotes}

\usepackage[ruled, linesnumberedhidden,vlined]{algorithm2e} 
\usetikzlibrary{matrix,decorations.markings}
\usepackage{amsmath,amssymb,dsfont,mathtools}
\usepackage{hyperref}
\hypersetup{
    colorlinks=true,       
    linkcolor=blue,          
    citecolor=olive,       
    filecolor=magenta,      
    urlcolor=blue          
}
\usepackage{aliascnt}

\usepackage{soul}

\usepackage[round,comma,sort,longnamesfirst]{natbib}

\DeclareMathOperator{\lcm}{lcm}

\begin{document}

\newcommand{\ca}{\mathcal{C}_{AL}}
\newcommand{\discn}{\mathbb{D}_{n}}
\title{Canonical Reduction Systems \\ in Artin-Tits groups of spherical type}
\date{ }
\author{Mar\'{i}a Cumplido, Juan Gonz\'alez-Meneses, Davide Perego}


\maketitle
\theoremstyle{plain}
\newtheorem{theorem}{Theorem}

\newaliascnt{lemma}{theorem}
\newtheorem{lemma}[lemma]{Lemma}
\aliascntresetthe{lemma}
\providecommand*{\lemmaautorefname}{Lemma}

\newaliascnt{proposition}{theorem}
\newtheorem{proposition}[proposition]{Proposition}
\aliascntresetthe{proposition}
\providecommand*{\propositionautorefname}{Proposition}

\newaliascnt{corollary}{theorem}
\newtheorem{corollary}[corollary]{Corollary}
\aliascntresetthe{corollary}
\providecommand*{\corollaryautorefname}{Corollary}

\newaliascnt{conjecture}{theorem}
\newtheorem{conjecture}[conjecture]{Conjecture}
\aliascntresetthe{conjecture}
\providecommand*{\conjectureautorefname}{Conjecture}

\theoremstyle{remark}

\newaliascnt{claim}{theorem}
\newaliascnt{remark}{theorem}
\newtheorem{claim}[claim]{Claim}
\newtheorem{remark}[remark]{Remark}
\newaliascnt{notation}{theorem}
\newtheorem{notation}[notation]{Notation}
\aliascntresetthe{notation}
\providecommand*{\notationautorefname}{Notation}

\aliascntresetthe{claim}
\providecommand*{\claimautorefname}{Claim}

\aliascntresetthe{remark}
\providecommand*{\remarkautorefname}{Remark}

\newtheorem*{claim*}{Claim}
\theoremstyle{definition}

\newaliascnt{definition}{theorem}
\newtheorem{definition}[definition]{Definition}
\aliascntresetthe{definition}
\providecommand*{\definitionautorefname}{Definition}

\newaliascnt{example}{theorem}
\newtheorem{example}[example]{Example}
\aliascntresetthe{example}
\providecommand*{\exampleautorefname}{Example}

\def\autorefspace{\hspace*{-0.5pt}}
\def\sectionautorefname{Section\autorefspace}
\def\subsectionautorefname{Section\autorefspace}
\def\subsubsectionautorefname{Section\autorefspace}
\def\figureautorefname{Figure\autorefspace}
\def\subfigureautorefname{Figure\autorefspace}
\def\tableautorefname{Table\autorefspace}
\def\equationautorefname{Equation\autorefspace}
\def\Itemautorefname{item\autorefspace}
\def\Hfootnoteautorefname{footnote\autorefspace}
\def\AMSautorefname{Equation\autorefspace}

\newcommand{\co}{\simeq_c}
\newcommand{\w}{\widetilde}
\newcommand{\po}{\preccurlyeq}
\newcommand{\dist}{\mathrm{d}}

\newcommand{\firstpartialcurve}[1]{%
  \draw[
    postaction={
      decorate,
      decoration={
        markings,
        mark=between positions 0 and #1 step 0.01 with {
          \draw[line width=4.1pt,white] (-\pgflinewidth,0) -- (\pgflinewidth,0);
        }
      }
    }
  ]
  (3,0) to[out=-90,in=90] (1,-2);
}

\newcommand{\secondpartialcurve}[2]{%
  \draw[
    postaction={
      decorate,
      decoration={
        markings,
        mark=between positions #1 and #2 step 0.01 with {
          \draw[line width=4.1pt,white] (-\pgflinewidth,0) -- (\pgflinewidth,0);
        }
      }
    }
  ]
  (1,2) to[out=-90,in=90] (3,0);
}

\newcommand{\thirdpartialcurve}[2]{%
  \draw[
    postaction={
      decorate,
      decoration={
        markings,
        mark=between positions #1 and #2 step 0.01 with {
          \draw[line width=3.5pt,white] (-\pgflinewidth,0) -- (\pgflinewidth,0);
        }
      }
    }
  ]
  (2.5,0) to[out=-90,in=90] (3.5,-2);
}

\def\Z{\mathbb Z} 
\def\Ker{{\rm Ker}} \def\R{\mathbb R} \def\GL{{\rm GL}}
\def\HH{\mathcal H} \def\C{\mathbb C} \def\P{\mathbb P}
\def\SSS{\mathfrak S} \def\BB{\mathcal B} 
\def\supp{{\rm supp}} \def\Id{{\rm Id}} \def\Im{{\rm Im}}
\def\MM{\mathcal M} \def\S{\mathbb S}
\newcommand{\bigveer}{\bigvee^\Lsh}
\newcommand{\wedger}{\wedge^\Lsh}
\newcommand{\veer}{\vee^\Lsh}
\def\diam{{\rm diam}}

\newcommand{\myref}[2]{\hyperref[#1]{#2~\ref*{#1}}}

\begin{abstract}
We introduce the canonical reduction system of an element in an Artin-Tits group of spherical type, which generalizes the similar notion for braids (and mapping classes) introduced by Birman, Lubotzky and McCarthy. We show its basic properties, which coincide with those satisfied in braid groups, and we provide an algorithm to compute it. We improve the algorithm in the case of braid groups, and discuss its complexity in this case. As a necessary result for obtaining the general algorithm, we prove that the centralizers of positive powers of an element form a periodic sequence and we show how to compute its period.

\medskip

{\footnotesize

 }

\end{abstract}

\section{Introduction}

The braid group on \( n \) strands, introduced by \cite{Artin1947}, is the mapping class group of the \( n \)-punctured disk. Algebraically, it has the following presentation:

$$
\mathcal{B}_n = \left\langle \sigma_1, \dots, \sigma_{n-1} \, \middle| \,
\begin{array}{lr}
\sigma_i \sigma_j = \sigma_j \sigma_i, & |i - j| > 1 \\
\sigma_i \sigma_j \sigma_i = \sigma_j \sigma_i \sigma_j, & |i - j| = 1
\end{array}
\right\rangle .
$$
Both the topological perspective—viewing the action of the braid group on the curve complex of the \( n \)-punctured disk, $\discn$,—and the combinatorial perspective—using Garside theory—have been key tools in understanding the structure of braids over the past few decades. From the topological point of view, the Nielsen-Thurston classification allows us to classify elements of the mapping class group of a surface into three distinct types: periodic, reducible (non periodic) and pseudo-Anosov. 

\medskip

A mapping class may preserve a family of (isotopy classes of) disjoint, simple closed curves in the surface, which we will call a system of curves. Cutting the surface along these curves, the restriction of~\(g\) to each connected component is a \emph{reduction} of~\(g\) into simpler mapping classes. A system of curves is \emph{adequate} for~\(g\) if such a reduction produces only periodic and pseudo-Anosov mapping classes. 
In the 1980s, \cite{BLM} proved that for every mapping class~\( g \), the set of adequate systems for \(g\) have a unique minimal element under inclusion. This minimal element is called the \emph{canonical reduction system} of~\(g\), denoted \(CRS(g)\). A mapping class~\(g\) is reducible, non-periodic, if and only if \(CRS(g)\neq \varnothing\). The canonical reduction system is an important tool for proving distinct algebraic results in mapping class groups and particularly in braid groups.

\medskip

In the case of braid groups, one can also use its Garside structure, an algebraic structure of the group, which allows to compute normal forms and solve the conjugacy problem, among other properties. Some generalizations of braid groups, such as Artin–Tits groups of spherical type, also admit a Garside structure, so any Garside-theoretic argument that does not rely on braid-specific properties will also be valid for these groups.

\medskip

Some of the topological objects which are defined for braid groups (seen as mapping classes) have an algebraic translation. Namely, an isotopy class of simple closed curves in the punctured disk corresponds to an \emph{irreducible parabolic subgroup}, and the complex of curves of the punctured disk $\mathcal{C}(\discn)$ is analogous to the complex of irreducible parabolic subgroups $\mathcal{C}(\mathcal{B}_n) $\citep{CGGW}. The action of a braid on a curve (as a mapping class) is equivalent to the action on the corresponding parabolic subgroup (by conjugation).  

\medskip

The translation of topological notions and arguments to algebraic ones, would allow to extend results from braid groups to Artin groups of spherical type. This paper aims to take a step in this direction: given an Artin-Tits group~$G$ of spherical type, and an element $\alpha\in G$, we introduce a completely algebraic definition of $CRS(\alpha)$, the canonical reduction system of~$\alpha$, which coincides with the classical definition if $G$ is a braid group.  We prove that this algebraic notion of $CRS(\alpha)$ satisfies similar properties as the classical one: it is preserved by powers and by multiplication by central elements, and it behaves as expected under conjugations.

Also, we provide an algorithm to compute $CRS(\alpha)$ (\autoref{alg_computing_CRS_Artin}). One of the key ingredients to obtain this algorithm is the fact, shown in this paper (\autoref{th_centralizers_are_periodic}), that the sequence of centralizers $Z(\gamma),Z(\gamma^2),Z(\gamma^3),\ldots$ for an element $\gamma\in G$ is periodic, and that one can compute its period.

Finally, we provide a better algorithm to compute $CRS(\alpha)$ in the particular case of braid groups (\autoref{alg_computing_CRS_braids}), and we discuss its complexity.

\medskip

In the particular case of braids, other tools can be used to detect reducing curves, like the theory of train-tracks \citep{BestvinaHandel1995}, and some algorithms have been implemented using that approach \citep{Hall}, but they do not compute the canonical reduction system, as far as we know. We apply and extend the results from \citep{Benardete1993}, the distinct (Garside-theoretic) solutions to the conjugacy problem in Artin-Tits groups of spherical type, and we uncover relations between reducing curves (or subgroups) and centralizers, to obtain our algorithms.

\medskip

The article is structured in the following manner: in Section~2 we give all the necessary background on the complex of curves, the complex of irreducible parabolic subgroups, Garside theory and mapping class groups. We also show that some of the results relating the complex of curves and Garside theory in braids can be extended to all Artin-Tits groups of spherical type. In Section~3 we define the canonical reduction system of any element $\alpha$, and prove its main properties. In Section~4 we give the necessary tools to obtain the algorithms to compute canonical reduction systems, we undertake a deep study of the centralizers of powers of any element $\alpha$, and we give the general algorithm to compute $CRS(\alpha)$ (\autoref{alg_computing_CRS_Artin}). Finally, in Section~5 we treat the case of braid groups, introducing \autoref{alg_computing_CRS_braids} and discussing its complexity.

\section{Definitions and background}

\subsection{Curve complex and canonical reduction system}

Given a surface $S$ (possibly with boundary and with punctures), we will consider isotopy classes of simple closed curves in the interior of $S$ that are non-degenerate, that is, they are not homotopic to a point, or a puncture, or a boundary component. We say that two such isotopy classes of curves are disjoint if we can find representatives in their respective classes that are disjoint. For simplicity, we will call each of these classes a \emph{curve}. The \emph{curve complex}~$\mathcal C(S)$ of the surface $S$ \citep{Harvey} is a flag complex having as vertices the curves of~$S$ with the following adjacency condition: two vertices share an edge if the corresponding curves are disjoint. Hence, a simplex of the curve complex of dimension~$d$ is a collection of $d+1$ mutually disjoint curves. It is known that one can take $d+1$ representatives, one for each curve, so that any two of them have empty intersection. For technical reasons, we will also admit the emptyset as a simplex of~$\mathcal C(S)$ of dimension $-1$. Notice that a mapping class permutes the vertices of~$\mathcal C(S)$, respecting adjacencies, so it determines an isomorphism of the curve complex.

\medskip

The Nielsen-Thurston classification arranges the elements of the mapping class group of a surface in three types, depending on the action they induce on the curve complex of the surface. An element~\( g \) is \emph{periodic} if the action induced by $g$ has finite order; it is \emph{reducible} if it fixes some simplex of the curve complex; and it is \emph{pseudo-Anosov} otherwise. Since reducible elements can be periodic, one usually classifies the mapping classes as periodic, reducible (non-periodic), and pseudo-Anosov. In this way one has a partition on the set of mapping classes of the given surface.
In the pseudo-Anosov case, there exist two transverse measured foliations on the surface for which the distance between leaves is scaled by factors~\( \lambda \) and~\( \frac{1}{\lambda} \), respectively, under the action of~\( g \) \citep{Thurston1988}. Moreover, the action induced by a pseudo-Anosov element on the curve complex does not have any finite orbit (that is, no curve is preserved by a nontrivial power of a pseudo-Anosov element).

\medskip

In the case of braid groups, which are torsion-free, the only elements that induce a trivial action on the complex of curves are the powers of \(\Delta^2\), where \(\Delta\) is the half-twist (or Garside element) \(\Delta=\sigma_1(\sigma_2\sigma_1)(\sigma_3\sigma_2\sigma_1)\cdots (\sigma_{n-1}\cdots \sigma_1)\). It is well known that the center of \(\mathcal{B}_n\) is \(Z(\mathcal{B}_n)=\langle\Delta^2\rangle\).  Hence, we have the following classification: a braid~$\alpha$ is periodic if it has a non-trivial central power; it is reducible if it fixes some simplex of the curve complex; and it is pseudo-Anosov if every vertex of the curve complex has an infinite orbit.

\medskip

A simplex of the curve complex $\mathcal C(S)$ will be called a \emph{multicurve} or a \emph{system of curves}. Given a  mapping class~\(g\) of a surface $S$, a \emph{reduction system} of $g$ is a multicurve which is preserved by~\(g\). Notice that $g$ is reducible (possibly being periodic too) if and only if it admits a nonempty reduction system. 

\medskip

Given a nonempty reduction system \(M\) of \(g\), we can take some $m>0$ so that each curve in~\(M\) is fixed by~\(g^m\). Then we can consider the connected components obtained by cutting the surface along~$M$, and the restriction of~\(g^m\) to each of these components. In this way, we \emph{reduce} (a power of) the original mapping class into several mapping classes. We say that \(M\) is an \emph{adequate} reduction system if every mapping class obtained by this decomposition is either periodic or pseudo-Anosov. According to~\cite{Thurston1988}, every reducible, non-periodic mapping class~\(g\) admits an adequate reduction system. If~\(g\) is either periodic or pseudo-Anosov, it is clear that the empty set is an adequate reduction system for~\(g\). Hence, every mapping class admits an adequate reduction system.

\medskip

In the 1980s, \cite{BLM} proved that for every mapping class~\( g \), the set of adequate reduction systems for~\(g\) has a unique minimal element under inclusion. This minimal element is called the \emph{canonical reduction system} of~\(g\), denoted \(CRS(g)\). It can also be seen as the intersection of all maximal (with respect to the inclusion) adequate reduction systems of~\(g\). A mapping class~\(g\) is reducible, non-periodic, if and only if \(CRS(g)\neq \varnothing\). 

\medskip

In the case of the braid group, the decomposition of the punctured disk along a multi\-curve yields a family of punctured disks. Hence, a reducible braid can be reduced into simpler braids (braids with less strands). The reduction along the canonical reduction system allows to show results concerning braids by first proving them for the periodic and pseudo-Anosov cases (see~\citealp{BLM}, \citealp{GM2003}, \citealp{GMWiest}).


\subsection{Complex of irreducible parabolic subgroups}

In addition to being a particular case of mapping class groups, braid groups are also a particular case of Artin-Tits groups. More precisely, they are Artin-Tits groups of spherical type, which are a kind of groups having very convenient algebraic properties (see next subsection on Garside theory). 

\medskip

\begin{figure}[h]
  \centering
  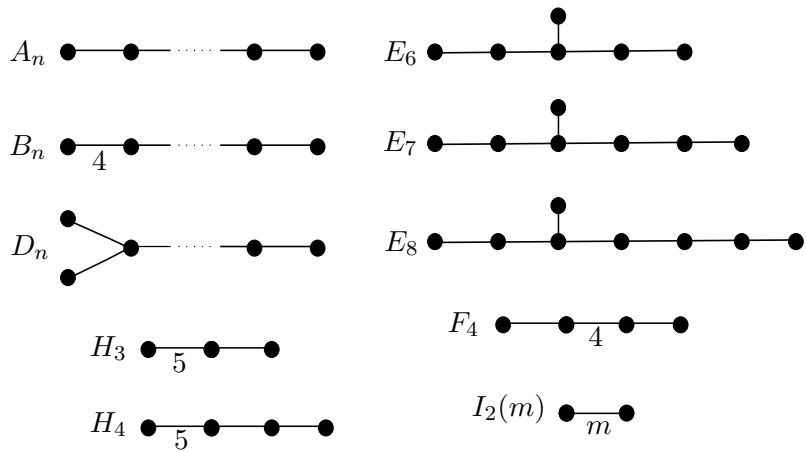
  \medskip
  \caption{The complete classification of irreducible Artin-Tits groups of spherical type.}\label{Dibujo:Spherical-type}
\end{figure}

Artin–Tits groups of spherical type, classified by their Dynkin diagrams (see \autoref{Dibujo:Spherical-type}), are generated by a finite set of standard generators. These generators correspond to the vertices of the Dynkin diagram; in the case of the braid group~$\mathcal{B}_{n}$, also denoted~$A_{n-1}$, they are $\sigma_1, \ldots, \sigma_{n-1}$. The defining relations of the group are determined by the edges of the diagram. If there is no edge between two vertices corresponding to generators $s$ and $t$, then $s$ and $t$ commute. If they are joined by an unlabeled edge, we have $sts = tst$. If they are joined by an edge labeled $m$, we have $stst\cdots = tsts\cdots$, where each side of the equality has $m$ letters.

Given a notion or property of braids defined via their interpretation as mapping classes, it is natural to ask whether it extends to all Artin–Tits groups of spherical type. One important example is the action of braid groups on the curve complex, which can be generalized to all Artin-Tits groups of spherical type, as shown in \citep{CGGW}. A curve (that is, an isotopy class of non-degenerate, simple closed curves in the punctured dis~$\discn$) corresponds to a subgroup of the braid group~$\mathcal{B}_n$, consisting of mapping classes supported in the punctured disk enclosed by the curve. This subgroup turns out to be an \emph{irreducible parabolic subgroup} of~$\mathcal{B}_n$. We now explain this notion.

Given~$G$, an Artin-Tits group of spherical type, a \emph{standard parabolic subgroup} of $G$ is a subgroup generated by a subset of the standard generators. Such a subgroup is \emph{irreducible} if its generators are adjacent in the Dynkin diagram, that is, if the full subgraph determined by those vertices is connected. In the case of the braid group, for which the Dynkin diagram is a linear graph, this just means that the generators are consecutive: \(\{\sigma_i, \sigma_{i+1}, \ldots, \sigma_j\}\). 

An (irreducible) \emph{parabolic subgroup} of $G$ is a conjugate of an (irreducible) standard parabolic subgroup. The trivial subgroup and the whole group $G$ are parabolic subgroups, but we will be interested in \emph{proper} parabolic subgroups, distinct from these two. In the case of braid groups, it is easy to see that there is a bijection between the (isotopy classes of non-degenerate, simple closed) curves of $\discn$ and the proper irreducible parabolic subgroups of $\mathcal{B}_n$. 

\medskip

Recall that the definition of the curve complex relied on the notion of \emph{adjacency} of curves. This can also be translated to algebraic terms, as shown in \citep{CGGW}. An irreducible parabolic subgroup $P$ of $G$ has a cyclic center, with a unique generator $z_P$ which is conjugate to a positive element. Two irreducible parabolic subgroups $P_1$ and $P_2$ are said to be \emph{adjacent} if $z_{P_1}$ and $z_{P_2}$ commute. In the case of braid groups, this happens if and only if the corresponding curves are disjoint, that is, if and only if they are adjacent in the curve complex \citep{CGGW}.

\medskip

Hence, every statement in the context of the braid group $\mathcal{B}_n$ which involves simple closed curves and disjointness, can be extended to all Artin-Tits groups of spherical type. On the contrary, if a statement involves the restriction of a braid (as a mapping class) to a subsurface of the punctured disk, then it cannot be easily generalized to other Artin-Tits groups. In particular, the curve complex can then be naturally generalized:

\begin{definition}[{\citealp{CGGW}}]
Let $G$ be an Artin-Tits group of spherical type. The complex of irreducible parabolic subgroups of $G$, denoted $\mathcal C(G)$, is a simplicial complex in which a simplex of dimension $d$ is a set $M=\{P_0,\ldots, P_d\}$ of proper irreducible parabolic subgroups such that $z_{P_i}$ commutes with $z_{P_j}$ for every $0\leq i, j\leq d$ (that is, the subgroups in $M$ are mutually adjacent).
\end{definition}

Notice that we will also consider $M=\varnothing$ as a simplex of dimension $-1$ of $\mathcal C(G)$.

\medskip

In the same way that braids act on the curve complex as mapping classes, every element of $G$ acts on the complex of irreducible parabolic subgroups by conjugation. The action is the same in the case of $\mathcal{B}_n$, so we have extended to all Artin-Tits groups of spherical type not only the notion of curve complex, but also the action of the elements of the group on it.

We finish this section with an important property.

\begin{theorem}[{\citealp[Theorem 2.2]{CGGW}}]\label{theorem_adjacency}
  Let $G$ be an Artin-Tits group of spherical type. Two irreducible parabolic subgroups $P_1$ and $P_2$ of $G$ are adjacent if and only if one of the following properties is satisfied: 
  \begin{itemize}
      \item $P_1\subset P_2$.
      \item $P_2\subset P_1$.
      \item $P_1\cap P_2=\{1\}$ and $\gamma_1\gamma_2=\gamma_2\gamma_1$ for every $\gamma_1\in P_1$ and $\gamma_2\in P_2$.
  \end{itemize} 
\end{theorem}

\subsection{Garside Theory and $\mathcal C(G)$}

It is well-known that an Artin-Tits group $G$ of spherical type has a Garside structure \citep{Brieskorn1972,Dehornoy1999}, which allows to define normal forms. Let~$G^+$ be the monoid of positive elements (the submonoid of $G$ generated by the standard generators). Given two elements $\alpha,\beta\in G$, we say that $\alpha$ is a prefix of $\beta$, and we write $\alpha\preccurlyeq \beta$, if $\alpha^{-1}\beta\in G^+$. We say that $\alpha$ is a suffix of $\beta$, and we write $\beta \succcurlyeq \alpha$, if $\beta \alpha^{-1}\in G^+$. The Garside structure is characterized by the following properties: 

\begin{itemize}

 \item The partial order $\preccurlyeq$, which is clearly invariant under left-multiplication, is also a lattice order. That is, for every $a,b\in G$ there exists a unique greatest common divisor $a \wedge b$ and a unique least common multiple $a \vee b$ with respect to~$\po$.
 
\item There is an element $\Delta$ (in this case, the least common multiple of all standard generators), called the Garside element of $G$, which satisfies $\Delta^{-1} {G^+} \Delta={G^+}$.

\item The positive prefixes of $\Delta$  coincide with the positive suffixes of $\Delta$. They are called \emph{simple elements}, and they are a finite set of generators of $G$.

\item $G^+$ is atomic: for every positive element, the lengths of all positive words in the standard generators (also called atoms) representing it have an upper bound. In this case, all representatives of a positive element have the same length, so the upper bound is just the length of any representative. 

\end{itemize}

By the symmetry of the relations in the standard presentation of $G$, one can see that the suffix order $\succcurlyeq$, which is invariant under right-multiplication, is also a lattice order. The greatest common divisor and least common multiple of $a$ and $b$ with respect to $\succcurlyeq$ are denoted $a\wedge^{\Lsh} b$ and $a\vee^{\Lsh}b$, respectively. 

We will say that a simple element is \emph{proper} if it is distinct from $1$ and $\Delta$. That is, if $1\prec s \prec \Delta$ or, equivalently, if $\Delta \succ s \succ 1$.

It is well known that, in every Garside group, a positive power of $\Delta$ is central. In the case of Artin-Tits groups of spherical type, either $\Delta$ or $\Delta^2$ is central.

\begin{definition}[\citealp{Elrifai1994}]
A decomposition $\beta=\Delta^k \beta_1\cdots \beta_N$ of an element $\beta\in G$ is the \emph{left normal form} of $\beta$ if 

\begin{itemize}

\item $k\in \mathbb{Z}$;

\item $\beta_i$ is a proper simple element for $i=1,\ldots , N$;

\item  $\beta_i \beta_{i+1}$ is in left normal form for $0<i < N$, that is, for every non-trivial positive prefix $a\po \beta_{i+1}$, the braid $\beta_ia$ is not simple.

\end{itemize}

\end{definition}

The above normal form, which can be computed in quadratic time with respect to the length of $\beta$, allows to solve the word problem in the groups $G$. The existence of the above normal form also shows that for every $\beta\in G$ there exists some $m\geq 0$ such that $\Delta^m \beta$ is positive. It suffices to take $m=\max\{-k,0\}$.

An alternative statement of the third property in the above definitions is: $\beta_i\beta_{i+1} \wedge \Delta = \beta_i$ for $i=1,\ldots,N-1$. Another alternative statement is: $\beta_i\cdots \beta_N \wedge \Delta= \beta_i$ for $i=1,\ldots,N-1$.

\medskip

The Garside structure in $G$ allows to solve the conjugacy problem by computing, for any given element, a finite invariant set of its conjugacy class. From the distinct possible invariant sets that have been proposed, in this article we will use the smallest one, the set of sliding circuits. We first need some definitions. 

\begin{definition}[\citealp{Gebhardt2010a}]
Let $\beta=\Delta^k \beta_1\cdots \beta_N$ be in left normal form. If $N>0$, we call $\iota(\beta):=\Delta^k \beta_1 \Delta^{-k}$ the \emph{initial factor} of $\beta$ and $\varphi(\beta):=\beta_N$ the \emph{final factor} of $\beta$. If $N=0$, we define $\iota(\beta):=1$ and $\varphi(\beta):=\Delta$. 
\end{definition}

It is not difficult to show that the initial and final factors of $\beta$ and $\beta^{-1}$ are closely related. Namely, $\varphi(\beta)\iota(\beta^{-1})=\Delta=\varphi(\beta^{-1})\iota(\beta)$. 

\begin{definition}[\citealp{Gebhardt2010a}]
Given $\beta\in G$, the \emph{preferred prefix} of~$\beta$ is $\mathfrak{p}(\beta):=\iota(\beta)\wedge \iota(\beta^{-1})$. The \emph{cyclic sliding} of $\beta$ is the conjugate of $\beta$ by its preferred prefix: $\mathfrak{s}(\beta):=\mathfrak{p}(\beta)^{-1}\beta\mathfrak{p}(\beta)$.
\end{definition}

Since $\varphi(\beta)\iota(\beta^{-1})=\Delta$, we can describe the preferred prefix $\mathfrak{p}(\beta)$ as the longest positive prefix of $\iota(\beta)$ such that $\varphi(\beta)\mathfrak{p}(\beta)$ is still simple.  

\begin{definition}[\citealp{Gebhardt2010a}]
We say that $\gamma$ belongs to a sliding circuit if $\mathfrak{s}^m(\gamma)=\gamma$ for some $m>0$. The \emph{set of sliding circuits} $SC(\beta)$ of  $\beta$  is the set of all conjugates of $\beta$ belonging to a sliding circuit.
\end{definition}

In the case of the braid group, we can embed the punctured disk~$\discn$ in the complex plane in such a way that the punctures lie on the real axis. We say that a curve in $\discn$ is \emph{standard} if it only intersects the real axis in two points. By a suitable isotopy, a standard curve can be represented by a geometric circle. Recall that a curve has to be non-degenerate, so it must enclose more than one and less than $n$ punctures. Notice that a curve is standard if and only if its corresponding irreducible parabolic subgroup is standard.

 We say that a multicurve $M$ is standard if all its curves are standard (we also declare $M=\varnothing$ as a standard multicurve). Given a curve $C$ (resp. a multicurve $M$), we will denote by $C^g$ (resp. $M^g$) its image under the action of $g$.  Notice that the braid $\Delta=\sigma_1(\sigma_2\sigma_1)(\sigma_3\sigma_2\sigma_1)\cdots (\sigma_{n-1}\cdots\sigma_1)$ has the same action on the curves of $\discn$ as a rotation of $\discn$ by an angle of $\pi$. Hence, if $C$ (resp. $M$) is standard then $C^{\Delta^k}$ (resp $M^{\Delta^k}$) is standard for every $k\in \mathbb Z$.

\begin{theorem}[{\citealp[Theorem~5.7]{Benardete1995}, \citealp[Theorem~3.8]{Lee2008}}]\label{theorem_normal_form_preserves} Let $M$ be a standard multicurve and let $\beta=\Delta^k \beta_1\cdots \beta_N$ be a braid in left normal form. If $M^\beta$ is standard, then $M^{\Delta^k \beta_1\cdots \beta_i}$ is standard for $i\in\{0,\dots, N\}$.
\end{theorem}

The following result is closely related to the previous one.
We say that a \emph{positive} braid~$\alpha$ is a \emph{standardizer} of a multicurve $M$ if $M^\alpha$ is standard. The set of standardizers of $M$ is denoted~$\operatorname{St}(M)$.

\begin{theorem}[{\citealp[Theorem~4.2]{Lee2008}}] \label{theorem_wedge_preserves}
Let $M$ be a multicurve of $\discn$. The set $\operatorname{St}(M)$ is nonempty and closed under $\wedge$ and $\vee$, hence it is a sublattice of $(\mathcal{B}_n^+,\preccurlyeq)$. Therefore, there is a unique minimal element in $\operatorname{St}(M)$ with respect to $\preccurlyeq$, called the \emph{minimal standardizer} of $M$.
\end{theorem}

The last two results can be extended to all Artin-Tits groups of spherical type. In this case we say that a \emph{positive} element $\alpha$ is a \emph{right standardizer} of a parabolic subgroup $P$ if $P^\alpha := \alpha^{-1}P\alpha$ is standard. The set of right standardizers of $P$ is denoted $\operatorname{St}(P)$. In the same way, we say that a positive element $\alpha$ is a \emph{left standardizer} of $P$ if $\alpha P \alpha^{-1}$ is standard. The set of left standardizers of $P$ is denoted $\operatorname{St}^{\Lsh}(P)$.

\begin{theorem}[{\citealp[Corollary 2]{Cumplido2019b}}]\label{th_standardizer_of_parabolic_is_lattice}
Let $P$ be a parabolic subgroup of an Artin-Tits group $G$ of spherical type. The set $\operatorname{St}(P)$ of right standardizers of $P$ is nonempty and closed under $\wedge$ and $\vee$, hence it is a sublattice of $(G^+,\preccurlyeq)$. Therefore, there is a minimal element of $\operatorname{St}(P)$ with respect to $\preccurlyeq$, called the \emph{minimal standardizer} of $P$.
\end{theorem}

We remark that a simplex $M$ of $\mathcal C(G)$ does not necessarily correspond to a parabolic subgroup, since the irreducible parabolic subgroups that form $M$ could be nested. In any case, we can extend the above result to all simplices. A \emph{right standardizer} of a simplex $M=\{P_0,\ldots,P_d\}$ is a positive element $\alpha$ such that $M^\alpha=\{P_0^{\alpha},\ldots,P_d^{\alpha}\}$ consists of standard irreducible parabolic subgroups.

\begin{corollary}\label{cor_standardizer_of_simplex_is_lattice}
Let $G$ be an Artin-Tits group of spherical type and let $M$ be a simplex of $\mathcal C(G)$. The set $\operatorname{St}(M)$ of right standardizers of $M$ is closed under $\wedge$ and $\vee$, hence it is a sublattice of $(G^+,\preccurlyeq)$. 
\end{corollary}

\begin{proof}
Let $M=\{P_0,\ldots,P_d\}$. The result follows from \autoref{th_standardizer_of_parabolic_is_lattice} since, by definition,
$$
 \operatorname{St}(M)=\bigcap_{i=0}^{d}{\operatorname{St}(P_i)}.
$$
\end{proof}

It is important to notice that, a priori, $\operatorname{St}(M)$ could be empty. We will show below that this is never the case. Because of the symmetry of the defining relations of~$G$, the set $\operatorname{St}^{\Lsh}(M)$ of \emph{left standardizers} of $M$ is closed under $\wedge^{\Lsh}$ and $\vee^{\Lsh}$, and hence forms a sublattice of $(G^+,\succcurlyeq)$. We begin with the case of simplices of~$\mathcal{C}(G)$ whose vertices are not related by inclusion. Notice that if~$P$ is a parabolic subgroup, a standardizer of~$P$ clearly exists: some conjugate of $P$ is standard by definition, and the conjugating element can be chosen positive after multiplication by a suitable central power of~$\Delta$. In contrast, for a simplex with no inclusion relations we obtain a collection of mutually commuting irreducible parabolic subgroups, which, a priori, need not be simultaneously standardized by a single conjugating element. We now show that this is indeed the case, and therefore such a simplex actually corresponds to a parabolic subgroup.

\begin{proposition}\label{prop_max_standardizer_nonempty}
Let $G$ be an Artin-Tits group of spherical type. Let $M=\{P_0,\ldots,P_d\}$ be a simplex of $\mathcal C(G)$ whose vertices do not satisfy any inclusion relation (that is, $P_i\not\subset P_j$ for every $i\neq j$). 
Then $\operatorname{St}(M)\neq \varnothing$. Furthermore, the minimal subgroup $P=\langle P_0,\ldots, P_d\rangle$ of $G$ containing $P_0,\ldots,P_d$ is a parabolic subgroup of $G$ whose irreducible components are $P_0,\ldots, P_d$, that is, $P$ is isomorphic to $P_0\times \cdots \times P_d$.
\end{proposition}

\begin{proof}
For $d=-1$ the result is trivial (\(\operatorname{St}(\varnothing)=G^+\)), and for $d=0$ it follows from \autoref{th_standardizer_of_parabolic_is_lattice}. Assume henceforth that $d>0$ and that the result holds for all smaller values of $d$.
 We can then \emph{standardize} the subgroups $P_0,\ldots,P_{d-1}$ using a common positive conjugating element $\beta$, hence we can consider the simplex $\{P_0^{\beta},\ldots,P_{d-1}^{\beta},P_d^{\beta}\}$, which will satisfy the result if and only if the original simplex satisfies it. Therefore, we can assume that the subgroups $P_0,\ldots, P_{d-1}$ are \emph{standard} irreducible parabolic subgroups.

Given an element $\delta\in G$, its pn-normal form is a decomposition $\delta=ab^{-1}$ such that $a,b$ are positive elements whose maximal common suffix is trivial, that is, $a\wedge^{\Lsh} b =1$. Recall that the lattice order $\wedge^{\Lsh}$ is invariant under right-multiplication, hence we have $ab^{-1}\wedge^{\Lsh} 1 = b^{-1}$. That is, $b^{-1}=\delta\wedge^{\Lsh} 1$.

We need to find an element $\gamma$ such that $P_i^{\gamma}$ is standard for every $i=0,\ldots,d$. Recall from \autoref{th_standardizer_of_parabolic_is_lattice} that $P_d$ admits a minimal standardizer. In~\cite[Theorem 3]{Cumplido2019b} it is shown that the minimal standardizer of $P_d$ is $b$, where $a b^{-1}$ is the pn-normal form of $z_{P_d}$. That is, $b^{-1}=z_{P_d}\wedge^{\Lsh} 1$. Let us see that we can take $\gamma=b$.

We know that $(P_d)^b$ is standard, because $b$ is the minimal standardizer of $P_d$. So we need to show that $(P_i)^b$ is also standard for every $i<d$. 

Let $L$ be a positive integer, big enough so that $z_{P_d}\Delta^{L}$ is positive and such that $\Delta^{L}$ is central in $G$. Then, since $b^{-1}=z_{P_d}\wedge^{\Lsh} 1$ and $\wedge^{\Lsh}$ is preserved under right-multiplication, we obtain $b^{-1}\Delta^L = z_{P_d}\Delta^L\wedge^{\Lsh} \Delta^L$. Hence $b^{-1}\Delta^L$ is positive, being the greatest common suffix of two positive elements.

Now let $i\in\{0,\ldots,d-1\}$. Since $P_d$ and $P_i$ are adjacent and they do not satisfy an inclusion relation, we obtain from \autoref{theorem_adjacency} that every element of $P_d$ commutes with every element of $P_i$. Hence, $z_{P_d}$ centralizes $P_i$. Since $\Delta^L$ is central, this implies that both $z_{P_d}\Delta^L$ and $\Delta^L$ centralize the standard parabolic subgroup~$P_i$. This means that $z_{P_d}\Delta^L$ and~$\Delta^L$ belong to~$\operatorname{St}^{\Lsh}(P_i)$. Since~$\operatorname{St}^{\Lsh}(P_i)$ is invariant under~$\wedge^{\Lsh}$, it follows that $\left(b^{-1}\Delta^L\right)P_i \left(b^{-1}\Delta^L\right)^{-1}$ is standard. But this latter element is no other than $b^{-1}P_ib$, as $\Delta^L$ is central.
Therefore, $b^{-1}P_i b$ is standard, as we wanted to show. That is, we can take $\gamma=b$ as an element which conjugates every $P_i$, for $i=0,\ldots,d$, to a standard parabolic subgroup.

It follows that conjugation by $\gamma$ sends the group $P=\langle P_0,\ldots, P_d\rangle$ to $Q=\langle Q_0,\ldots, Q_d\rangle$ where $Q_i=(P_i)^{\gamma}$ is a proper, irreducible standard parabolic subgroup, for $i=0,\ldots,d$. Since every~$Q_i$ is generated by a subset of the standard generators of~$G$, the subgroup~$Q$ is a standard parabolic subgroup, hence~$P$ is a parabolic subgroup. Moreover, since every~$Q_i$ is irreducible, the standard generators of~$Q_i$ are adjacent in the Dynkin graph of $G$. Finally, since every element of $Q_i$ commutes with every element of~$Q_j$ for $i\neq j$, and $Q_i\cap Q_j=1$, no standard generator of~$Q_i$ coincides or is adjacent to a standard generator of~$Q_j$. It follows that $Q_0,\ldots,Q_d$ are the irreducible components of~$Q$, and that~$Q$ is isomorphic to $Q_0\times \cdots \times Q_d$. Hence (conjugating by~$\gamma^{-1}$), the irreducible components of $P$ are $P_0, \ldots, P_d$ and~$P$ is isomorphic to $P_0\times \cdots \times P_d$.
\end{proof}

Let us now prove the existence of a standardizer (hence a minimal standardizer) in the general case. It should be noted that alternative proofs of the following two results was previously provided by \cite[Proposition~3.8, Proposition~3.13]{Ragosta}.

\begin{proposition}\label{prop_standardizer_nonempty}
Let $G$ be an Artin-Tits group of spherical type. For every  simplex~$M$ of~$\mathcal C(G)$, the set $\operatorname{St}(M)$ of right standardizers of~$M$ is nonempty, hence there exists a minimal right standardizer of~$M$. The set $\operatorname{St}^{\Lsh}(M)$ of left standardizers of $M$ is also nonempty, hence there exists a minimal left standardizer of~$M$.
\end{proposition}

\begin{proof}
Let $M=\{P_0,\ldots,P_d\}$ a simplex. By symmetry of the standard relations in $G$, we just need to show the result for $\operatorname{St}(M)$. If this set is nonempty, \autoref{cor_standardizer_of_simplex_is_lattice} will assure the existence of a minimal right standardizer of $M$. 

Let us then show that $\operatorname{St}(M)$ is nonempty. This is trivial if $d=-1$ and already known if $d=0$, so we will assume that $d>0$ and that the result holds for smaller values of $d$.

Let $Q_1,\ldots,Q_r$ be the maximal subgroups in $M$, with respect to inclusion. Then $\{Q_1,\ldots,Q_r\}$ satisfies the hypothesis of \autoref{prop_max_standardizer_nonempty}, so there exists some $\alpha\in \operatorname{St}(\{Q_1,\ldots,Q_r\})$. We can then consider $M^\alpha$ and assume that all maximal subgroups in $M$ are standard.

For every $i=1,\ldots,r$, the group $Q_i$ is an Artin-Tits group of spherical type, whose Dynkin diagram is just the maximal subgraph of the graph of $G$ determined by the standard generators of $Q_i$. Let $M_i$ be the subset of $M$ consisting of those subgroups strictly contained in $Q_i$. We know by \citep[Theorem~0.2]{Godelle} that every parabolic subgroup contained in $Q_i$ is a parabolic subgroup of the group $Q_i$. Then $M_i$ is a simplex of $\mathcal C(Q_i)$, with less than $d+1$ vertices. By induction hypothesis, there is some $\gamma_i\in Q_i$ which standardizes $M_i$. Since the subgroups $Q_1,\ldots,Q_r$ mutually commute, it follows that the element $\gamma:=\gamma_1\cdots \gamma_r\in G$ standardizes $M$. Hence $\operatorname{St}(M)$ is nonempty. 
\end{proof}

An important consequence is the following:

\begin{corollary}\label{cor_dimension_C(G)}
Let $G$ be an Artin-Tits group of spherical type. Let $\Sigma$ be the set of standard generators of $G$. Then $\dim(\mathcal C(G))=\#(\Sigma)-2$.
\end{corollary}

\begin{proof}
The fact that $\mathcal C(G)$ is finite-dimensional was already observed in \citep[Corollary~3.9]{Ragosta}. Moreover, the equality $\dim(\mathcal C(G)) = \#(\Sigma) - 2$ follows immediately from \citep[Proposition~3.13]{Ragosta}, where maximal simplices are characterized. We provide below an alternative proof.

Let $\Sigma=\{s_1,\ldots,s_n\}$. If $n=1$, then $G$ is cyclic, so there are no proper standard parabolic subgroups, and hence there are no proper irreducible parabolic subgroups. Therefore, $\mathcal C(G)$ is empty and we have $\dim(\mathcal C(G))=-1=\#(\Sigma)-2$ in this case.

Let us then suppose that $n>1$ and that the result holds for smaller values of $n$. 

The set
$$
  \{\{s_1\},\{s_1,s_2\},\ldots,\{s_1,\cdots,s_{n-1}\}\}
$$
is a simplex of $\mathcal C(G)$ of dimension $n-2$. Hence $\dim(\mathcal C(G))\geq \#(\Sigma)-2$. Let us show the converse inequality.

Let $M$ be a simplex of $\mathcal C(G)$. By \autoref{cor_standardizer_of_simplex_is_lattice}, $M$ can be standardized, so we can assume that all vertices of $M$ are standard.

Let $\{P_1,\ldots,P_r\}$ be the set of maximal subgroups (by inclusion) in $M$, and let $\Sigma_i$ be the set of standard generators of $P_i$. Then $\Sigma_i\cap \Sigma_j=\varnothing$ for every $i\neq j$, Therefore $\Sigma_1\sqcup \cdots \sqcup \Sigma_r \subsetneq  \Sigma$. The equality cannot be achieved, as no element in $\Sigma_i$ can be adjacent to an element in $\Sigma_j$ and, if $r=1$, $P_1$ cannot be the whole group.

Now let $M_i$ be the set of subgroups of $M$ strictly contained in $P_i$. Again by \citep[Theorem~0.2]{Godelle}, $M_i$ is a simplex of $\mathcal C(P_i)$, hence $\#(M_i)< \#(\Sigma_i)$ by induction hypothesis. We then have $\#(M_i\cup \{P_i\})\leq \#(\Sigma_i)$ and, therefore, 
$$
\#(M)=\sum_{i=1}^{r}\#(M_i\cup \{P_i\}) \leq \sum_{i=1}^{r}\#(\Sigma_i) < \#(\Sigma).
$$
That is, $\dim(M)\leq \#(\Sigma)-2$ for every simplex $M$. Hence $\dim(\mathcal C(G))\leq \#(\Sigma)-2$.
\end{proof}

We end this section with two consequences of \autoref{cor_standardizer_of_simplex_is_lattice}, which generalize the corresponding results in $\mathcal{B}_n$.

\begin{corollary}\label{cor_BGN_for_parabolic}
Let $\beta=\Delta^k \beta_1\cdots \beta_N\in G$ be written in left normal form. If $M$ is a standard simplex of $\mathcal C(G)$ and $M^\beta$ is standard, then $M^{\Delta^k \beta_1\cdots \beta_i}$ is standard for every $i=0,\ldots, N$. 
\end{corollary}

\begin{proof}
For $i=0$ the result holds, since conjugation by $\Delta$ permutes the standard generators of~$G$, hence $M^{\Delta^k}$ is standard (each of its vertices is generated by standard generators of $G$).

Suppose that $i>0$ and that $M_i:=M^{\Delta^{k}\beta_1\cdots \beta_{i-1}}$ is standard. We know that $$\beta_i=(\beta_i\cdots \beta_N)\wedge \Delta.$$
Now $(M_i)^{\beta_i\cdots \beta_N}= M^{\Delta^k \beta_1\cdots \beta_N}=M^\beta$ is standard and $(M_i)^{\Delta}$ is also standard. Hence, since~$\operatorname{St}(M_i)$ is closed under $\wedge$, it follows that $(M_i)^{\beta_i}=(M_i)^{(\beta_i\cdots \beta_N)\wedge \Delta}$ is standard. That is, $M^{\Delta^k \beta_1\cdots \beta_i}$ is standard, as we wanted to show.
\end{proof}

As a further consequence, we have the following corollary:

\begin{corollary}\label{cor_cyclic_sliding_preserves_standard}
Let $G$ be an Artin-Tits group of spherical type and let $\beta\in G$. If $M$ is a standard simplex of $\mathcal C(G)$ preserved by $\beta$, then $M^{\mathfrak{p}(\beta)}$ is a standard simplex preserved by $\mathfrak{s}(\beta)$.
\end{corollary}

\begin{proof}
If $\beta$ is a power of $\Delta$ then $\mathfrak{p}(\beta)=1$ and the result is trivial. Let us then assume that the left normal form of $\beta$ is $\Delta^k \beta_1\cdots \beta_N$ for some $N\geq 1$. We have that $M^\beta=M$, so $M^\beta$ is standard. By \autoref{cor_BGN_for_parabolic}, $M^{\Delta^k\beta_1}$ is standard. Hence, $M^{\iota(\beta)}=\left(M^{\Delta^k \beta_1}\right)^{\Delta^{-k}}$ is standard. 

In the same way, as $M^{\beta^{-1}}=M$, we also have that $M^{\iota(\beta^{-1})}$ is standard. Therefore, by \autoref{cor_standardizer_of_simplex_is_lattice}, $\mathfrak p(\beta)=\iota(\beta)\wedge \iota(\beta^{-1})\in \operatorname{St}(M)$. This shows that $M^{\mathfrak{p}(\beta)}$ is a standard simplex.

The fact that this standard simplex is preserved by $\mathfrak s(\beta)$ is obvious, since
$$
  \left(M^{\mathfrak{p}(\beta)}\right)^{\mathfrak s(\beta)}=
  \left(M^{\mathfrak{p}(\beta)}\right)^{\mathfrak{p}(\beta)^{-1} \beta \mathfrak{p}(\beta)}=
  M^{\beta \mathfrak{p}(\beta)}=
    \left(M^\beta\right)^{\mathfrak{p}(\beta)}=
    M^{\mathfrak{p}(\beta)}.
$$
\end{proof}

\section{Canonical reduction systems in Artin-Tits groups}

Let us now try to generalize the notion of canonical reduction system to all Artin-Tits groups of spherical type. As we said before, the definition of the canonical reduction system of a braid (or a mapping class) is the following, given by~\cite{BLM}:

\begin{definition}\label{def_CRS_1}
 Given a braid $\alpha\in \mathcal{B}_n$, the canonical reduction system of $\alpha$, denoted $CRS(\alpha)$, is the minimal adequate reduction system of $\alpha$, with respect to inclusion.    
\end{definition}

Since the notion of \emph{adequate reduction system} involves restrictions of $\alpha$ to subsurfaces of $\discn$, it cannot be easily generalized to other Artin-Tits groups. Moreover, it is not clear how one could compute the canonical reduction system of a braid using the above definition.

Fortunately, \cite{BLM} have a criterion to determine if a given curve belongs to $CRS(\alpha)$, where $\alpha$ is a mapping class. We will state the definitions and results only for braids:

\begin{definition}[\citealp{BLM}]\label{def_CRS_2}
 Given a braid $\alpha\in \mathcal{B}_n$, a \emph{reduction curve} of $\alpha$ is a curve that belongs to a reduction system of $\alpha$. Equivalently, a reduction curve of $\alpha$ is a curve $C$ whose orbit under $\alpha$ is a simplex of $\mathcal C(\discn)$. A reduction curve $C$ of $\alpha$ is \emph{essential} if it is disjoint from any curve $C'$ whose orbit under $\alpha$ is finite.
\end{definition}

\begin{theorem}[{\citealp[Lemma 2.6]{BLM}}]
  Given $\alpha\in \mathcal{B}_n$, its canonical reduction system $CRS(\alpha)$ is the set of essential reduction curves of $\alpha$.  
\end{theorem}

Notice that the latter characterization of a canonical reduction system can be easily extended to all Artin-Tits groups of spherical type, as it only involves curves fixed by (some power of) $\alpha$, and disjointness of curves.  As the group $G$ acts on $\mathcal C(G)$ by conjugation, we can define:

\begin{definition}
Let $G$ be an Artin-Tits group of spherical type. Given $\alpha\in G$, a \emph{reduction simplex} of $\alpha$ is a simplex $M$ of $\mathcal C(G)$ which is invariant under the action of $\alpha$ (i.e. conjugation by $\alpha$ permutes the subgroups in $M$). A \emph{reduction subgroup} of $\alpha$ is a vertex of a reduction simplex of $\alpha$. Equivalently, a reduction subgroup of $\alpha$ is a proper irreducible parabolic subgroup of $G$ whose orbit under $\alpha$ is a simplex of $\mathcal C(G)$. A reduction subgroup $P$ of $\alpha$ is \emph{essential} if it is adjacent to any other irreducible parabolic subgroup whose orbit under $\alpha$ is finite.
\end{definition}

\begin{definition}
 Let $G$ be an Artin-Tits group of spherical type, and let $\alpha\in G$. The \emph{canonical reduction system} of $\alpha$, $CRS(\alpha)$, is the set of essential reduction subgroups of $\alpha$.   
\end{definition}

To prove that $P \in CRS(\alpha)$, it suffices to show that the orbit of $P$ under $\alpha$ is finite, and that $P$ is adjacent to any parabolic subgroup $Q$ with finite orbit under $\alpha$. The latter condition ensures that the orbit of $P$ forms a simplex.   

\medskip

It is clear that, in the case of braid groups, the two given definitions of $CRS(\alpha)$ are equivalent, using the natural bijection between non-degenerate simple closed curves and proper irreducible parabolic subgroups.

\begin{proposition}
Let $G$ be an Artin-Tits group of spherical type, and let $\alpha\in G$. The $CRS(\alpha)$ is a simplex of $\mathcal C(G)$. 
\end{proposition}

\begin{proof}
A reduction subgroup of $\alpha$ has a finite orbit under the action of $\alpha$. Hence, if $P_1$ and $P_2$ are essential reduction subgroups of $\alpha$, by the fact that $P_1$ is essential and $P_2$ has a finite orbit under $\alpha$, then $P_1$ and $P_2$ are adjacent. It follows that if $CRS(\alpha)$ is nonempty, it is formed by mutually adjacent, proper irreducible parabolic subgroups. 

Moreover, if $CRS(\alpha)$ is nonempty, every finite subset of $CRS(\alpha)$ is a simplex of $\mathcal C(G)$.  By \autoref{cor_dimension_C(G)}, the cardinality of every simplex of $\mathcal C(G)$ is bounded above by the number of standard generators of $G$. Therefore, $CRS(\alpha)$ is finite, so it is a simplex of $\mathcal C(G)$.
\end{proof}

Let us see that, in the context of Artin-Tits groups of spherical type, the canonical reduction system satisfies some of the properties which are satisfied in $\mathcal{B}_n$ (for the analogous properties in $\mathcal{B}_n$, see~\citealp[Lemma 2.6]{BLM}).

\begin{proposition}\label{prop_CRS_conjug_power_center}
Let $G$ be an Artin-Tits group of spherical type, and let $\alpha\in G$. One has:
\begin{enumerate}
  \item $CRS(\alpha)=CRS(\alpha^m)$ for every $m\in \mathbb Z\setminus\{0\}$.
  \item $CRS\left(\alpha^\beta\right)=CRS(\alpha)^{\beta}$ for every $\beta\in G$.
  \item $CRS(z \alpha)=CRS(\alpha)$ for every $z\in Z(G)$, the center of $G$.
\end{enumerate}
\end{proposition}

\begin{proof}\phantom{x}

\begin{enumerate}
  \item We will denote by $P^{\langle \delta\rangle}$ the orbit of a subgroup $P$ under the action of an element $\delta$. 
  
  Let $P\in CRS(\alpha)$, that is, let $P$ be an essential reduction subgroup of $\alpha$, and let $m\neq 0$. We have $P^{\langle\alpha^m \rangle} \subset P^{\langle\alpha\rangle}$, and hence, since $P^{\langle\alpha\rangle}$ is a simplex of $\mathcal C(G)$, $P^{\langle\alpha^m \rangle}$ is also a simplex of $\mathcal C(G)$. Therefore, $P$ is a reduction subgroup of $\alpha^m$. 

  Let $Q$ be a vertex of $\mathcal C(G)$ such that $Q^{\langle\alpha^m\rangle}$ is finite. This implies that $Q^{\langle\alpha\rangle}$ is also finite (at most $m$ times bigger), hence $P$ and $Q$ are adjacent, as $P$ is an essential subgroup of~$\alpha$. Therefore, $P$ is an essential subgroup of $\alpha^m$. In other words, $P\in CRS(\alpha^m)$.

  Conversely, let $P\in CRS(\alpha^m)$, that is, let $P$ be an essential reduction subgroup of $\alpha^m$ for some $m\neq 0$. Then $P^{\langle\alpha^m \rangle}$ is a simplex. In particular, $P^{\langle\alpha^m \rangle}$ is finite, so $P^{\langle\alpha \rangle}$ is also finite. We must first show that it is a simplex.

  If $P^{\langle\alpha \rangle}$ is not a simplex, there are $i\neq j$ such that $P^{\alpha^i}$ is not adjacent to $P^{\alpha^j}$. Conjugating both subgroups by $\alpha^{-i}$, we obtain that $P$ and $P^{\alpha^{j-i}}$ are not adjacent. But the orbit of $P^{\alpha^{j-i}}$ under $\alpha^m$ is finite (it is contained in $P^{\langle\alpha\rangle}$), so this contradicts that $P$ is essential for $\alpha^m$. Therefore, $P^{\langle\alpha \rangle}$ is a simplex, so $P$ is a reduction subgroup of $\alpha$.

  Now let $Q$ be a vertex of $\mathcal C(G)$ such that $Q^{\langle\alpha\rangle}$ is finite. This 
  implies that $Q^{\langle\alpha^m\rangle}$ is finite, and hence, as $P$ is an essential reduction subgroup of $\alpha^m$, $P$ and $Q$ are adjacent. Therefore, $P$ is an essential reduction subgroup of $\alpha$. In other words, $P\in CRS(\alpha)$.

 \item  The conjugacy by $\beta$ preserves adjacency, so it induces an isomorphism between $\mathcal{C}(G)$ and~$\mathcal{C}(G)^\beta$. In terms of actions, we have $P^\alpha=Q$ if and only if $\left(P^\beta\right)^{\alpha^\beta}=Q^\beta$.
 The definition of $CRS(\alpha)$ only depends on the action of $\alpha$ and adjacency conditions. Then, $P$ is an essential reduction subgroup of $\alpha$ if and only if $P^\beta$ is an essential reduction subgroup of~$\alpha^\beta$.

 \item Since $z$ is central, the actions induced by $z\alpha$ and $\alpha$ on $\mathcal{C}(G)$ coincide. The canonical reduction system of an element is determined by its action on $\mathcal{C}(G)$, and therefore the result follows.
\end{enumerate}
\end{proof}

\section{Detecting vertices of the canonical reduction system}

\subsection{Restricting the search to a finite set}

We know that the canonical reduction system of an element $\alpha\in G$ is a simplex of $\mathcal C(G)$. In general, this complex has infinitely many vertices. We now describe how to compute a finite set of vertices of~$\mathcal C(G)$ that contains $CRS(\alpha)$. In fact, the set we will compute is a simplex of~$\mathcal C(G)$, so it will have at most~$\#(\Sigma)-1$ elements, where $\Sigma$ is the set of standard generators of~$G$, by \autoref{cor_dimension_C(G)}.

We will use a similar approach to that undertaken in \citep{GMWiest11} for braids. For every element $\gamma\in G$, let
$$
  R_{\gamma}=\{M \mid  M \text{ is a standard reduction simplex of }\gamma\}.
$$
Recall that the empty set is a standard reduction simplex, hence $R_{\gamma}$ is always nonempty, as it has $M=\varnothing$ as an element. Now we define $d_{\gamma}=\max\{\dim(M) \mid  M\in R_{\gamma}\}$. Notice that $-1\leq d_{\gamma}\leq \#(\Sigma)-2$.

\begin{proposition}\label{prop_max_simplex_contains_CRS}
Let $G$ be an Artin-Tits group of spherical type, and let $\alpha\in G$. Let $\gamma\in SC(\alpha)$ such that $d_{\gamma}$ is maximal, and let $M$ be a standard reduction simplex of $\gamma$ with dimension $d_{\gamma}$. Then $CRS(\gamma)$ is a subset of $M$, and $CRS(\alpha)$ is a subset of $M^{c}$, where $c$ is any conjugating element such that $\gamma^c=\alpha$.
\end{proposition}

\begin{proof}
Recall from \autoref{prop_CRS_conjug_power_center} that, if $\gamma=\alpha^{\beta}$ for some conjugating element $\beta$, then $CRS(\gamma)=CRS(\alpha)^{\beta}$. Hence, if $CRS(\alpha)=\varnothing$ we will have $CRS(\gamma)=\varnothing$ for every $\gamma\in SC(\alpha)$. In this case, the result is trivially true.

Let us then assume that $CRS(\alpha)\neq \varnothing$. By \autoref{prop_standardizer_nonempty}, the standardizer of $CRS(\alpha)$ is nonempty, so we can consider an element $\beta_1$ such that $CRS(\alpha)^{\beta_1}=CRS(\alpha^{\beta_1})$ is standard.  Now we will apply iterated cyclic sliding to $\alpha^{\beta_1}$ until we obtain an element $\alpha^{\beta_1\beta_2}\in SC(\alpha)$, where $\beta_2$ is the product of all preferred prefixes by which one conjugates. By \autoref{cor_cyclic_sliding_preserves_standard}, $CRS(\alpha^{\beta_1})^{\beta_2}=CRS(\alpha^{\beta_1\beta_2})$ is standard. Therefore, there exists some element in $SC(\alpha)$ whose (nonempty) canonical reduction system is standard. 

Let us consider $\gamma \in SC(\alpha)$ such that $d_{\gamma}$ is maximal, and let $M\in R_{\gamma}$ such that $\dim(M)=d_{\gamma}$. By the previous argument, $d_{\gamma}\geq 0$, hence~$M$ is nonempty.

Suppose that $CRS(\gamma)$ is not a subset of $M$. Then the set $M_+=CRS(\gamma)\cup M$ is bigger than~$M$. Moreover, $M_+$~is a simplex of $\mathcal C(G)$, as every two subgroups $P$ and $Q$ in $M_+$ are adjacent: if~$P$ and~$Q$ belong to $CRS(\gamma)$, they are adjacent because $CRS(\gamma)$ is a simplex; if they both belong to $M$, they are adjacent because~$M$ is a simplex; if $P\in CRS(\gamma)$ and $Q\in M\setminus CRS(\gamma)$, they are adjacent because $P$ is an essential reduction subgroup of $\gamma$ and the orbit of~$Q$ under~$\gamma$ is finite (as it belongs to a reduction simplex of~$\gamma$). Therefore, $M_+$ is a simplex of~$\mathcal C(G)$.

Moreover, $M_+$ is a reduction simplex of $\gamma$, as 
$$
(M_+)^{\gamma}=CRS(\gamma)^{\gamma}\cup M^{\gamma} = CRS(\gamma)\cup M=M_+.
$$
But, a priori, $M_+$ is not necessarily standard.

Let $c_1\in G^+$ be any right standardizer of $M^+$, which is known to exist thanks to \autoref{prop_standardizer_nonempty}. Then $(M_+)^{c_1}$ is a standard reduction simplex of $\gamma^{c_1}$. But now $\gamma^{c_1}$ does not necessarily belong to $SC(\alpha)$.

By \autoref{cor_cyclic_sliding_preserves_standard}, if we apply iterated cyclic sliding to $\gamma^{c_1}$ until we obtain an element $\widetilde{\gamma}=\gamma^{c_1c_2}\in SC(\alpha)$, the simplex $(M_+)^{c_1c_2}$ will be a standard reduction simplex of $\widetilde \gamma$. This implies that $d_{\widetilde\gamma}\geq \dim((M_+)^{c_1c_2})=\dim(M_+)>\dim(M)=d_{\gamma}$, with $\widetilde\gamma\in SC(\alpha)$. This contradicts the maximality of $d_\gamma$. Therefore, $CRS(\gamma)$ must be a subset of $M$.

Finally, if $c$ is any conjugating element such that $\gamma^c=\alpha$, then $CRS(\gamma)\subset M$ is equivalent to $CRS(\alpha)=CRS(\gamma)^c \subset M^c$.
\end{proof}

We can improve the above result very easily.

\begin{corollary}\label{cor_max_simplex_contains_CRS}
Let $G$ be an Artin-Tits group of spherical type, and let $\alpha\in G$. Let $\gamma\in SC(\alpha)$ such that $d_{\gamma}$ is maximal, let $M_1,\ldots,M_r$ be the standard reduction simplices of $\gamma$ with dimension~$d_{\gamma}$, and let $M=M_1\cap \cdots \cap M_r$. Then $CRS(\gamma)$ is a subset of~$M$, and $CRS(\alpha)$ is a subset of~$M^{c}$, where~$c$ is any conjugating element such that $\gamma^c=\alpha$.
\end{corollary}

\begin{proof}
 By \autoref{prop_max_simplex_contains_CRS}, $CRS(\gamma)$ is contained in $M_i$ for $i=1,\ldots, r$. Hence $CRS(\gamma)\subset M$. And this is equivalent to $CRS(\alpha)=CRS(\gamma)^c\subset M^c$.   
\end{proof}

We can restrict even more the set $M$, by applying the properties of canonical reduction systems.

\begin{proposition}\label{prop_CRS_preserved_by_centralizer}
 Let $\alpha\in G$ and let $M$ be a simplex of $\mathcal C(G)$ which contains $CRS(\alpha)$. Let $Z(\alpha)=\langle z_1,\ldots,z_m\rangle$ be the centralizer of $\alpha$, and let 
 $$
  M_0=M\cap M^{z_1} \cap \cdots \cap M^{z_m}.
 $$ 
 Then $CRS(\alpha)\subset M_0$.
\end{proposition}

\begin{proof}
 For every $z\in Z(\alpha)$, we have $$CRS(\alpha)^z=CRS(\alpha^z)=CRS(\alpha).$$ Hence, since $CRS(\alpha)\subset M$, we have $CRS(\alpha)=CRS(\alpha)^z\subset M^z$. This implies that $$CRS(\alpha)\subset M\cap M^{z_1}\cap \cdots \cap M^{z_m}=M_0.$$
\end{proof}

After \autoref{cor_max_simplex_contains_CRS} and \autoref{prop_CRS_preserved_by_centralizer}, we obtain a procedure to compute a simplex $M_0$ containing $CRS(\alpha)$, for every $\alpha\in G$. Moreover, that simplex will be invariant under the action of any element in $Z(\alpha)$. The procedure goes as follows:

\begin{enumerate}

\item  Compute the set $SC(\alpha)$ and, for every $\gamma\in SC(\alpha)$, a conjugating element $c_{\gamma}$ such that $\gamma^{c_\gamma}=\alpha$.

\item For every $\gamma\in SC(\alpha)$, compute the standard reduction simplices of $\gamma$ and determine $d_{\gamma}$.

\item Let $\gamma\in SC(\alpha)$ be such that $d_{\gamma}$ is maximal, and let $M$ be the intersection of the maximal standard reduction simplices of $\gamma$. 

\item Compute a set of generators $z_1,\ldots,z_m$ of $Z(\alpha)$, and let $M_0=M\cap M^{z_1}\cap\cdots \cap M^{z_m}$. Then~$(M_0)^{c_\gamma}$ contains $CRS(\alpha)$.
\end{enumerate}

Step~1 is a routine computation \citep{Gebhardt2010a,Gebhardt2010}. For Step~2, we need to compute the standard reduction simplices of an element $\gamma$. Once this is done, we can select the maximal ones and determine their intersection. Computing the centralizer of an element is also standard; see \citep{Franco2003b}, whose results apply to sets of sliding circuits, and also \citep{GMWiest,GMValladares}. The final assertion in Step~4 follows from \autoref{cor_max_simplex_contains_CRS} and \autoref{prop_CRS_preserved_by_centralizer}.

Let us then see how to compute the standard reduction simplices of an element.

\subsection{Computing standard reduction simplices}

Let $G$ be an Artin-Tits group of spherical type and let $\alpha\in G$. We will show a procedure to compute the standard reduction simplices of $\alpha$, based on results in \citep{Benardete1993} and \citep{CGGW}. 

\medskip

Notice that if $G$ is an Artin-Tits group of spherical type which is not irreducible, every irreducible parabolic subgroup will be contained in one component, and all its conjugates will remain in the same component, so we can restrict our attention to irreducible Artin-Tits groups of spherical type. For the rest of this section, $G$ will be irreducible.

\medskip

We can enumerate all irreducible standard parabolic subgroups in~$G$. The Dynkin graph of~$G$ is either linear, or becomes linear after removing a single vertex together with its unique adjacent edge We denote by $a_1,\ldots,a_r$ the vertices of the linear graph, and~$\widetilde{a}$ the extra vertex, if it exists. The \emph{rank} of~$G$ will be $n=\#(\Sigma)$, so that $n=r$ when the graph of $G$ is linear, and $n=r+1$ otherwise.

If the graph of $G$ is linear ($r=n$), its irreducible standard parabolic subgroups are $A_{[i,j]}=\langle a_i,a_{i+1},\ldots,a_j\rangle$ for $1\leq i\leq j\leq n$ except for $(i,j)=(1,n)$. Hence, the number of irreducible parabolic subgroups is $\frac{n(n+1)}{2}-1 = \frac{n^2+n-2}{2}$. 

If the graph of $G$ is not linear, we have either $G=D_n$ for $n\geq 4$ or $G\in\{E_6,E_7,E_8\}$. In this case its irreducible parabolic subgroups are $A_{[i,j]}=\langle a_i,a_{i+1},\ldots,a_j\rangle$ for $1\leq i\leq j\leq r$, also $A_{[i,j]}^+=\langle a_i,a_{i+1},\ldots,a_j,\widetilde{a}\rangle$ for $1\leq i\leq j \leq r$ except for $(i,j)=(1,r)$, provided that $\widetilde{a}$ is adjacent to one of the remaining generators, and $A_{\varnothing}^+=\langle \widetilde{a}\rangle$. Hence, in this case, the number of irreducible parabolic subgroups is smaller than $r(r+1)$, that is, smaller than $(n-1)n$.

The following table shows the cardinality of the set $\mathcal P_{IS}$ of irreducible standard parabolic subgroups for each irreducible Artin-Tits group.
$$
\begin{array}{|c||c|c|c|c|c|c|c|c|c|c|}
  \hline
    G &  A_n & B_n & D_n & E_6 & E_7 & E_8 & F_4 & H_3 & H_4 & I_{2}(p)
   \\ \hline 
    \#(\mathcal P_{IS}) & \frac{n^2+n-2}{2} & \frac{n^2+n-2}{2} &  \frac{n^2+3n-8}{2} & 24 & 33 & 43 & 9 & 5 & 9 & 2
   \\ \hline  
\end{array}
$$

Notice that if two irreducible standard parabolic subgroups $P_1=\langle s_1\cdots s_r\rangle$ and $P_2=\langle t_1\cdots t_m\rangle$ are conjugate (where each $s_i$ and $t_j$ are standard generators of $G$), then they are isomorphic. This implies that the subgraphs of the graph of $G$ they determine are isomorphic \citep{Paris2004}. In particular, we have $r=m$. For this reason, we will treat the standard parabolic subgroups of $G$ attending to its number of standard generators, that we will call its \emph{rank}.

We can see that, if the graph of $G$ is linear, the number of irreducible parabolic subgroups of rank $d$ is $n-d+1$. If the graph of $G$ is not linear, this number is at most $n-d+3$. In any case this number is at most $n$.

\medskip

For a given element $\alpha\in G$, we will first see how to compute the pairs of irreducible standard parabolic subgroups $(P,Q)$ such that $P^\alpha=Q$. 

Recall that, given an irreducible parabolic~$P$, the element~$z_P$ is the generator of its center which is conjugate to a positive element. Moreover, $P$~is standard if and only if $z_P$ is positive: Indeed, if $ab^{-1}$ is the pn-normal form of~$z_P$, then it is shown in \citep{Cumplido2019b} that $b$ is its minimal right standardizer. $P$~is standard if and only if its minimal standardizer~$b$ is trivial, and this happens if and only if $ab^{-1}$ is positive, that is, if~$z_P$ is positive.

We will use the following result:

\begin{theorem}[{\citealp[Lemma 8]{Cumplido2019b}}]\label{th_conjugate_parabolics_zP_zQ}
If $P$ and $Q$ are parabolic subgroups of $G$, and $\alpha\in G$, then $P^{\alpha}=Q$ if and only if $(z_P)^{\alpha}=z_Q$.
\end{theorem}

Thanks to this result, we just need to conjugate a single element in order to check whether a subgroup is conjugated into another one by an element~$\alpha$. Let us start with the case in which~$\alpha$ is a simple element:

\begin{lemma}\label{lem_standard_to_standard_by_simple}
Let $s$ be a simple element, and let $P$ be a standard irreducible parabolic subgroup. Then $P^s$ is standard if and only if $s\preccurlyeq z_P s$. In that case, $P^s=Q$ and $(z_P)^s=z_Q$, where $Q$ is the subgroup generated by the standard generators which are prefixes of $(z_P)^s$.
\end{lemma}

\begin{proof}
We know that $P^s$ is a parabolic subgroup, being the conjugate of a standard parabolic subgroup. If we denote $P^s=Q$, \autoref{th_conjugate_parabolics_zP_zQ} tells us that $(z_P)^s=z_Q$. Now $Q$ is standard if and only if $z_Q$ is positive, but $z_Q=s^{-1}z_Ps$, where $z_Ps$ is a positive element, so this happens if and only if $s\preccurlyeq z_Ps$.

Assume that $Q$ is standard. In order to obtain the standard generators of $Q$, we just need to recall that, since $Q$ is standard, $z_Q$ is either the Garside element $\Delta_Q$ (which is the least common multiple of the standard generators of $Q$) or $(\Delta_Q)^2$. Since the product $\Delta_Q \Delta_Q$ is in left normal form as written, the standard generators which are prefixes of $(\Delta_Q)^2$ are precisely the standard generators which are prefixes of $\Delta_Q$, and these are precisely the standard generators of $Q$.
\end{proof}

Notice that checking whether $P^s$ is standard is a very fast computation, as $z_Ps$ is a positive element with at most three simple factors, and we just need to compute its left normal form. Since $z_P$ is supposed to be already in normal form, this takes time $O(n^2)$, and check whether $s$ is a prefix of the first factor (also $O(n^2)$).

\medskip

Now we can treat the general case.

\begin{lemma}\label{lem_standard_to_standard}
Let $P$ be a standard irreducible parabolic subgroup, and let $\alpha=\Delta^k \alpha_1\cdots \alpha_N\in G$ be written in left normal form. Then $P^\alpha$ is standard if and only if, for $i=1,\ldots, N$, we have
$$
    \alpha_i\preccurlyeq z_{P_{i-1}}\alpha_i,
$$
where $z_{P_0}=(z_P)^{\Delta^k}$ and $z_{P_i}=(z_{P_{i-1}})^{\alpha_i}$ for $i=1,\ldots,N$.

In that case, $P^{\alpha}=Q$ and $(z_P)^\alpha=z_Q$, where $z_Q=z_{P_N}$ and $Q$ is the subgroup generated by the standard generators which are prefixes of $z_{P_N}$.
\end{lemma}

\begin{proof}
By definition, $P^\alpha$ is a parabolic subgroup, which we will call $Q$. By \autoref{cor_BGN_for_parabolic}, if $P^\alpha$ is standard, then $P^{\Delta^k\alpha_1\cdots \alpha_i}$ is standard for $i=0,\ldots, N$.  The converse is trivially true. If we denote $P_i=P^{\Delta^k\alpha_1\cdots \alpha_i}$, we have that $Q=P_N$ is standard if and only if $P_0,P_1,\ldots,P_N$ are all standard subgroups.

The subgroup $P_0=P^{\Delta^k}$ is always standard, since powers of $\Delta$ preserve the standardness of subgroups. Now notice that $P_i=(P_{i-1})^{\alpha_i}$, for $i=1,\ldots,N$, and that $z_{P_0}=(z_P)^{\Delta^k}$ and $z_{P_i}=(z_{P_{i-1}})^{\alpha_i}$ for $i=1,\ldots,N$. Since each $\alpha_i$ is a simple element, we can apply \autoref{lem_standard_to_standard_by_simple} to conclude that, assuming $P_{i-1}$ is standard, $P_i$ is standard if and only if $\alpha_i\preccurlyeq z_{P_{i-1}}\alpha_i$ for $i=1,\ldots, N$, as we wanted to show. 

If this is the case, $z_Q=z_{P_N}$ is a positive element (either $\Delta_Q$ or $(\Delta_Q)^2)$, and the standard generators which are prefixes of $z_{P_N}$ are precisely the standard generators of $Q$.
\end{proof}

We remark that checking whether $P^\alpha$ is standard takes time $O(Nn^2)$, as it performs $N$~times the computation of \autoref{lem_standard_to_standard_by_simple}. 

\medskip

Now we can easily obtain a procedure to compute the standard reduction simplices of $\alpha$. We need to compute the orbit under $\alpha$ of each standard irreducible parabolic subgroup $P$. Such an orbit is retained provided all its elements are standard and mutually adjacent.

\begin{algorithm}[h!]\caption{For computing the standard reduction simplices of an element in an Artin-Tits group $G$ of spherical type.}\label{algorithm_standard_simplices}\label{alg_standard_reduction_simplices}

\SetKwInOut{Input}{Input}\SetKwInOut{Output}{Output}

\BlankLine

\Input{An element $\alpha=\Delta^k\alpha_1\cdots \alpha_N\in G$ in left normal form.}
\Output{The set of nonempty standard reduction simplices of $\alpha$.}
\BlankLine
\BlankLine

Enumerate the irreducible standard parabolic subgroups $P_1,\ldots, P_r$ of $G$;

Set $U=\{P_1,\ldots,P_r\}$ and $V_1=\varnothing$;

\While{$U\neq \varnothing$}{
Extract $P$ from $U$ and set $\mathcal S=\{P\}$;

Compute $z_P$ and set $P_0=P$ and $z_{P_0}=z_P$;

\While{$\mathcal S\neq \varnothing$}{
 \If {$(P_0)^{\alpha}$ is standard (\autoref{lem_standard_to_standard})}{
  Set $z_Q=(z_{P_0})^{\alpha}$ and $Q=(P_0)^\alpha$;
  
  \If{$Q=P$}{
    Set $V_1=V_1\cup \{\mathcal S\}$;
    
    Set $\mathcal S=\varnothing$;
  }
  \Else{
  \If{$z_Qz_R=z_Rz_Q$ for every $R\in \mathcal S$}
  {
   Set $\mathcal S=\mathcal S\cup\{Q\}$; 

   Set $P_0=Q$ and $z_{P_0}=z_Q$;
  }
  \Else{
   Set $\mathcal S=\varnothing$;
  }
  }
 }
  \Else{
   Set $\mathcal S=\varnothing$;
  }
 }
}
Denote $\mathcal S_1,\ldots,\mathcal S_m$ the elements in $V_1$;

Set $V_2=\varnothing$;

\For{$i=1,\ldots,m-1$}{
 \For{$j=i+1,\ldots,m$}{
   \If{$z_Pz_Q=z_Qz_P$ for every $P\in \mathcal S_i$ and $Q\in \mathcal S_j$}{
   Set $V_2=V_2\cup \{\{i,j\}\}$;
  }
 }
}
Compute the set $CS$ of complete subgraphs of the graph with vertices $\{1,\ldots,m\}$ and edges $V_2$;

\Return{$\{\{\mathcal S_{i_1}\cup \cdots \cup \mathcal S_{i_l}\};\ \{i_1,\ldots,i_l\}\in CS\}$}

\end{algorithm}

In \autoref{algorithm_standard_simplices}, we take irreducible standard parabolic subgroups from the set $U$ which contains all of them, then we compute their orbits under $\alpha$ as long as the obtained elements are standard and disjoint to the previous ones. If this is the case, the orbit will arrive to the original subgroup and the orbit $\mathcal S$ will be stored. Otherwise, all standard subgroups contained in the orbit will be discarded. 

At the end of the first part of the algorithm, the set $V_1$ will contain all standard simplices $\mathcal S_1,\ldots,\mathcal S_m$ of $\mathcal C(G)$ which are orbits under the action of $\alpha$. In the second part of the algorithm, we check, for every $i<j$, whether $\mathcal S_i\cup \mathcal S_j$ is a simplex (that is, whether every two subgroups of $\mathcal S_i\cup \mathcal S_j$ are adjacent). We store the pairs $\{i,j\}$ for which $\mathcal S_i\cup \mathcal S_j$ is a simplex in the set $V_2$. Hence, the graph $\Gamma$ with set of vertices $\{1,\ldots,m\}$ and set of edges $V_2$ is the graph of adjacencies of the orbits $\mathcal S_1,\ldots,\mathcal S_m$. 

As every standard reduction simplex of $\alpha$ is precisely the union of orbits under $\alpha$, it follows that the set of standard reduction simplices of $\alpha$ is the set of complete subgraphs of $\Gamma$. Hence, the algorithm computes these complete subgraphs (one can use, for instance, the Bron-Kerbosch algorithm), and outputs the sets of subgroups corresponding to those complete subgraphs.

We observe that the ingredients involved in the complexity of \autoref{algorithm_standard_simplices} are:
\begin{itemize}

 \item The number of irreducible parabolic subgroups of $G$, which is $O(n^2)$. 

 \item Checking whether $(P_0)^\alpha$ is standard has complexity $O(Nn^2)$.

 \item Checking whether $z_Qz_P=z_Pz_Q$ takes time $O(n^2)$, as it computes and compares two normal forms of positive elements of length at most $4$.

 \item The total number of comparisons $z_Pz_Q=z_Qz_P$, in the whole algorithm (including the first and the second part) is $O(l^2)$, where $l$ is the number of standard irreducible parabolic subgroups in $G$. Hence it is $O(n^4)$. 

 \item Computing the complete subgraphs of $V_2$ has complexity $O(3^{n/3})$ \citep{TTT}.
\end{itemize}

By all these arguments, the total complexity of the algorithm is $O(Nn^6+3^{n/3})$.

\subsection{Discarding vertices which are not essential}

Given an element $\alpha\in G$, we have shown a procedure to obtain a simplex $M_0$ of $\mathcal C(G)$ containing $CRS(\alpha)$, which is invariant under the action of any element in $Z(\alpha)$. The remaining part of the algorithm to compute $CRS(\alpha)$ would be to detect which vertices of $M_0$ are essential, or, equivalently, to discard those which are not essential.

We will be able to do this efficiently in the case of braid groups, as we will see in \autoref{sec_braids}. But we can also give an algorithm that works in every Artin-Tits group of spherical type, thanks to the following result, involving centralizers of powers of elements.

\begin{theorem}\label{th_standard_CRS_and_centralizers}
Let $G$ be an Artin-Tits group of spherical type and let $\alpha\in G$. Let $\gamma$ be a conjugate of $\alpha$ such that $CRS(\gamma)$ is standard. Let $P$ be a standard reduction subgroup of $\gamma$. Then $P\in CRS(\gamma)$ if and only if $P^z$ is standard for every $z\in \bigcup_{m\geq 1} Z(\gamma^m)$.
\end{theorem}

\begin{proof}
Suppose that $P\in CRS(\gamma)$. Given $z\in \bigcup_{m\geq 1} Z(\gamma^m)$, let $r\geq 1$ be such that $z\in Z(g^r)$. Then, by \autoref{prop_CRS_conjug_power_center},
$$
  P^z\in CRS(\gamma)^z=CRS(\gamma^r)^z=CRS((\gamma^r)^z)=CRS(\gamma^r)=CRS(\gamma).
$$
Hence, as $CRS(\gamma)$ is standard, $P^z$ is standard.

Conversely, suppose that $P\notin CRS(\gamma)$. This implies that there is an irreducible parabolic subgroup $Q$ whose orbit under $\gamma$ is finite and is not adjacent to $P$. Let $r$ be the size of the orbit of $Q$ under $\gamma$. Then $Q^{\gamma^r}=Q$, hence $(z_Q)^{\gamma^r}=z_Q$. That is, $z_Q\in Z(\gamma^r)$.

Now, as $P$ is not adjacent to $Q$, we have $z_Pz_Q\neq z_Qz_P$. It is shown in \cite[Lemma 4.6]{CGGW} that this happens if and only if $(z_P)^l(z_Q)^m\neq (z_Q)^m(z_P)^l$ for every $l,m\neq 0$. Therefore, $(z_P)^{(z_Q)^m}\neq z_P$ for every $m>0$, and this implies that the sequence $\left\{(z_P)^{(z_Q)^m}\right\}_{m>0}$ does not have repeated elements, so it is infinite. Since the number of standard irreducible parabolic subgroups is finite, this implies that $P^{(z_Q)^m}$ is not standard for some $m>0$, where $(z_Q)^m\in Z(\gamma^r)$.
\end{proof}

In order to obtain an algorithm from the previous result, we need a finiteness condition, since in the statement we are considering an infinite union $\bigcup_{m\geq 1}{Z(\gamma^m)}$. The condition is given by the following result. We are deeply grateful to Ivan Marin, who suggested using linear representations, pointed out the importance of roots of unity, and provided a first proof of the result.

\begin{theorem}\label{th_centralizers_are_periodic}
Let $G$ be an Artin-Tits group of spherical type. For every $\gamma\in G$, the sequence $Z(\gamma),Z(\gamma^2),Z(\gamma^3),\ldots $ has a finite number of elements. Moreover, it is a periodic sequence, and the period can be computed from $\gamma$.
\end{theorem}

\begin{proof}
We will use the generalized Krammer representation of Artin-Tits groups of spherical type, given in \citep{CohenWales2002} and in \citep{Digne2003}. It is a faithful representation $\rho\colon G \rightarrow GL(W)$, where $W$ is a vector space over $\mathbb Q(t,q)$ of dimension $m$, the number of reflections in the Coxeter group determined by $G$ (alternatively, $m$ is the length of the Garside element $\Delta\in G$).

One can specialize $\rho$ by choosing particular complex values for $t$ and $q$, so that $\rho$ is still faithful. Hence, we can suppose that $\rho\colon G \rightarrow GL(\mathbb C^{m})$.

Let $\gamma \in G$, and denote $\mathbf{M}=\rho(\gamma)$. Since $\rho$ is faithful, the elements in $Z(\gamma)\subset G$ are in one-to-one correspondence (via $\rho$) with the elements in $Z(\mathbf{M})\cap \rho(G)\subset GL(\mathbb C^m)$. In the same way, the elements in $Z(\gamma^r)\subset G$ are in one-to-one correspondence with $$Z(\mathbf{M}^r)\cap \rho(G)\subset GL(\mathbb C^m).$$ Therefore, in order to prove the theorem, we just need to prove that the sequence $$Z(\mathbf{M}), Z(\mathbf{M}^2), Z(\mathbf{M}^3),\ldots $$ is finite (and, moreover, that it is periodic). The period of $\{Z(\gamma^r)\}_{r>0}$ will be a divisor of the period~$p$ of $\{Z(\mathbf{M}^r)\}_{r>0}$, hence the former can be obtained, once~$p$ is computed, by computing the centralizers $Z(\gamma),\ldots,Z(\gamma^p)$ and checking which elements are repeated. Therefore, our goal is to compute~$p$.

Now notice that the centralizer of a matrix is a vector subspace of $GL(\mathbb C^m)$, so it is a finite-dimensional $\mathbb C$-vector space. A classical result by Frobenius \cite[Page 111]{Jacobson} gives an explicit formula for this dimension in terms of the degrees of the invariant factors of the matrix. More precisely, if the degrees of the invariant factors of $\mathbf{M}$ are $d_1>d_2>\cdots > d_l$, then $\dim(Z(\mathbf{M}))=d_1+3d_2+\cdots+(2l-1)d_l$. It follows that the dimension of $Z(\mathbf{M})$ depends on the multiplicities of the eigenvalues of $\mathbf{M}$, and on the sizes of the Jordan blocks corresponding to each eigenvalue. But it does not depend on the actual eigenvalues of $\mathbf{M}$.

It will be useful to describe the above dimension in terms of the sizes of the Jordan blocks of $\mathbf{M}$. First, notice that we can rewrite the above sum as follows: 
$$
\dim(Z(\mathbf{M}))=\sum_{i=1}^{l}\sum_{j=1}^{l}{\min(d_i,d_j)}.
$$
Now let $\lambda_1,\ldots,\lambda_u$ be the distinct eigenvalues of $\mathbf{M}$, and let $n_{\lambda_h,1},\ldots,n_{\lambda_h,l}$ be the sizes of the Jordan blocks corresponding to $\lambda_h$, in descending order, where $n_{\lambda_h,j}=0$ if $j$ is bigger than the number of blocks corresponding to $\lambda_h$. Then, for $i=1,\ldots,l$ we have $d_i=n_{\lambda_1,i}+\cdots+n_{\lambda_u,i}$. By construction, if $d_i<d_j$ then $n_{\lambda_h,i}<n_{\lambda_h,j}$ for every eigenvalue $\lambda_h$. Hence $\min(d_i,d_j)=\min(n_{\lambda_1,i},n_{\lambda_1,j})+\cdots+\min(n_{\lambda_u,i},n_{\lambda_u,j})$. We can then separate in the above sum the contribution of each eigenvalue, and we obtain:
$$
\dim(Z(\mathbf{M}))=\sum_{h=1}^{u}\left(\sum_{i=1}^{l}\sum_{j=1}^{l}{\min(n_{\lambda_h,i},n_{\lambda_h,j})}\right).
$$
In the above decomposition, we may omit the zero summands by letting~$l_h$ denote the number of Jordan blocks for the eigenvalue $\lambda_h$, so that we obtain:
\begin{equation}\label{formula:dim_of_centralizer}
\dim(Z(\mathbf{M}))=\sum_{h=1}^{u}\left(\sum_{i=1}^{l_h}\sum_{j=1}^{l_h}{\min(n_{\lambda_h,i},n_{\lambda_h,j})}\right).
\end{equation}

We will also use the following result: If~$J(\lambda)$ is a Jordan block of size~$L$ corresponding to an eigenvalue $\lambda\neq 0$, then for every $r>0$, the matrix~$J(\lambda)^r$ is conjugate to the matrix~$J(\lambda^r)$. Indeed, one can easily show by induction that the matrix~$J(\lambda)^r$ has~$\lambda^r$ in every entry of the diagonal, and~$r\lambda^{r-1}$ in every entry of the super-diagonal. It follows that $J(\lambda)^r$ has~$\lambda^r$ as its only eigenvalue, and that $\operatorname{rank}(J(\lambda)^r-\lambda^rI)=L-1$. So the eigenspace of $J(\lambda)^r$ corresponding to~$\lambda^r$ has dimension 1. Therefore, the Jordan normal form of~$J(\lambda)^r$ has a single block of length~$L$, so~$J(\lambda)^r$ is conjugate to~$J(\lambda^r)$.

Let us go back to our sequence $\{Z(\mathbf{M}^r)\}_{r>0}$ of centralizers. To show that the sequence is periodic (hence finite), one can proceed as follows. Let $d=\max\{\dim(Z(\mathbf{M}^r))\}_{r>0}$, and let $p$ be the smallest positive integer such that $d=\dim(Z(\mathbf{M}^p))$. Since $Z(\mathbf{M}^p)\subset Z(\mathbf{M}^{rp})$ for every $r>0$, and the dimension of $Z(\mathbf{M}^p)$ is maximal, we have $Z(\mathbf{M}^p) = Z(\mathbf{M}^{rp})$ for every $r>0$. Moreover, for every $r>0$ we have $Z(\mathbf{M}^r)\subset Z(\mathbf{M}^{rp})=Z(\mathbf{M}^p)$, so every centralizer in the sequence is contained in $Z(\mathbf{M}^p)$.

It follows that, for every $r>0$, if $P\in Z(\mathbf{M}^r)$ then $P$ commutes with~$\mathbf{M}^r$ and with $\mathbf{M}^p$, so $P\in Z(\mathbf{M}^{r+p})$. Conversely, if $P\in Z(\mathbf{M}^{r+p})$ then $P$ commutes with $\mathbf{M}^{r+p}$ and also with $\mathbf{M}^{-p}$ (since $Z(\mathbf{M}^p)=Z(\mathbf{M}^{-p})$), hence $P\in Z(\mathbf{M}^r)$. Therefore, $Z(\mathbf{M}^{r})=Z(\mathbf{M}^{r+p})$ for every $r>0$, and we have shown that the sequence is periodic of period~$p$.

But this argument does not allow us to compute the period~$p$, as we do not know a priori which is the maximal dimension~$d$ in the sequence. Hence, we will use a different approach in order to compute the period.

First, notice that the matrix $\mathbf{M}=\rho(\gamma)$ is invertible, as $\mathbf{M}^{-1}=\rho(\gamma^{-1})$. Hence, all eigenvalues of $\mathbf{M}$ are nonzero. Notice also that, without loss of generality, after a suitable conjugation we can assume that $\mathbf{M}$ is in Jordan normal form: $\mathbf{M}=\mathbf{M}(\lambda_1)\oplus \cdots \oplus \mathbf{M}(\lambda_u)$, where $\lambda_1,\ldots,\lambda_u$ are the distinct (nonzero) eigenvalues of $\mathbf{M}$, and each $\mathbf{M}(\lambda_i)=J_1(\lambda_i)\oplus\cdots \oplus J_{d_i}(\lambda_i)$ is the direct sum of the $d_i$ Jordan blocks corresponding to $\lambda_i$, for $i=1,\ldots,u$.

For every $r>0$, we have $\mathbf{M}^r=\mathbf{M}(\lambda_1)^r\oplus \cdots \oplus \mathbf{M}(\lambda_u)^r$, where $$\mathbf{M}(\lambda_i)^r=J_1(\lambda_i)^r\oplus\cdots \oplus J_{d_i}(\lambda_i)^r,$$ for $i=1,\ldots, u$. Since $\lambda_i\neq 0$, we know that each $J_h(\lambda_i)^r$ is conjugate to $J_h(\lambda_i^r)$. Hence $\mathbf{M}(\lambda_i)^r$ is conjugate to $\mathbf{M}(\lambda_i^r)$ for $i=1,\ldots,u$. That is, $\mathbf{M}^r$ is conjugate to $\mathbf{M}(\lambda_1^r)\oplus \cdots \oplus \mathbf{M}(\lambda_u^r)$. Since the dimension of a centralizer only depends on the multiplicities of the eigenvalues, and the number and sizes of the Jordan blocks, it follows that, if $\lambda_1^r,\ldots,\lambda_u^r$ are all distinct, then $\dim(Z(\mathbf{M}))=\dim(Z(\mathbf{M}^r))$.

But the important observation is that $\lambda_1^r,\ldots,\lambda_u^r$ are not necessarily distinct. We have $\lambda_i^r=\lambda_j^r$ if and only if  $\left(\frac{\lambda_i}{\lambda_j}\right)^r=1$, so we just need to check the pairs of eigenvalues $\{\lambda_i,\lambda_j\}$ such that $\frac{\lambda_i}{\lambda_j}$ is a root of unity. Let $r_{ij}$ be the index of $\frac{\lambda_i}{\lambda_j}$ as a root of unity, if that is the case, and let $p$ be the least common multiple of all $r_{ij}$. We will see that $\dim(Z(\mathbf{M}^p))$ is maximal in the sequence $\{\dim(Z(\mathbf{M}^r))\}_{r>0}$, so the sequence is periodic of period $p$.

For every $r>0$, let $\mathcal P_r=\{(h_1,h_2):\ h_1\neq h_2 \text{ and } \lambda_{h_1}^r=\lambda_{h_2}^r\}$ be the set of ordered pairs of distinct eigenvalues which get identified when raised to the $r$-th power. Recall the formula~(\ref{formula:dim_of_centralizer}) for $\dim(Z(\mathbf{M}))$. If $\lambda_{h_1}^r=\lambda_{h_2}^r$, in the formula which determines $\dim(Z(\mathbf{M}^r))$ we have the same summands (as the blocks in the Jordan normal form of $\mathbf{M}^r$ have the same size as those in the Jordan normal form of $\mathbf{M}$), but we need to add new summands of the form $\min(m_1,m_2)$ corresponding to blocks coming from $\lambda_{h_1}$ and from $\lambda_{h_2}$, which now correspond to blocks of the same eigenvalue. Then we have
$$
   \dim(Z(\mathbf{M}^r))=\sum_{h=1}^{u}\left(\sum_{i=1}^{l_h}\sum_{j=1}^{l_h}{\min(n_{\lambda_h,i},n_{\lambda_h,j})}\right)
   +\sum_{(h_1,h_2)\in \mathcal P_r} \left(\sum_{i=1}^{l_{h_1}}\sum_{j=1}^{l_{h_2}} {\min(n_{\lambda_{h_1},i},n_{\lambda_{h_2},j})}\right).
$$
Finally, we just need to notice that in the formula for $\dim(Z(\mathbf{M}^p))$ we will see all summands which may appear in the formula for any other power, as all eigenvalues which can eventually be identified will be identified in $\mathbf{M}^p$. Therefore, $\dim(Z(\mathbf{M}^p))$ is maximal, as we wanted to show.

The number $p$ (the period of the sequence) can then be computed by checking which fractions~$\frac{\lambda_i}{\lambda_j}$ are roots of unity, computing their indices, and taking the least common multiple of all indices. By the above formula for the dimension, the number $p$ is the smallest possible period of the sequence $\{Z(\mathbf{M}^r)\}_{r>0}$. And from it one can obtain the smallest possible period of the sequence of centralizers of $\{Z(\gamma^r)\}_{r>0}$ (which is a divisor of $p$).
\end{proof}

The period of the sequence $Z(\gamma),Z(\gamma^2),Z(\gamma^3),\ldots$ depends a priori on the element $\gamma\in G$, as it depends on the eigenvalues of $\rho(\gamma)$. But we conjecture that there must be a common upper bound for all possible periods in an Artin-Tits group of spherical type.

\begin{conjecture}\label{conj_bound_for_period}
Let $G$ be an Artin-Tits group of spherical type. There is a number $n(G)$ such that, for every $\gamma\in G$, the period $p$ of the sequence $Z(\gamma),Z(\gamma^2),Z(\gamma^3),\ldots$ is at most $n(G)$.
\end{conjecture} 

Finding the bound $n(G)$ would provide an algorithm for finding $d$ with no need to compute the eigenvalues of $\rho(\gamma)$, as in the proof of \autoref{th_centralizers_are_periodic}. One would just need to compute the centralizers $Z(\gamma),Z(\gamma^2)\ldots, Z(\gamma^{n(G)})$.

In any case, the above results provide an algorithm to compute $CRS(\alpha)$ for an element $\alpha$ in an Artin-Tits group of spherical type. 

\bigskip 

\begin{algorithm}[H]\caption{For computing the canonical reduction system of an element in an Artin-Tits group of spherical type.}\label{algorithm_main_Artin}\label{alg_computing_CRS_Artin}

\SetKwInOut{Input}{Input}\SetKwInOut{Output}{Output}

\BlankLine

\Input{An element $\alpha\in G$, where $G$ is an Artin-Tits group of spherical type.}
\Output{$CRS(\alpha)=\mathcal{C}^c$ as a pair $(\mathcal C,c)$, where $\mathcal{C}$ is a standard simplex and $c\in  G$.}
\BlankLine
\BlankLine

Compute $SC(\alpha)$, and a conjugating element $c_\gamma$ such that $\gamma^{c_{\gamma}}=\alpha$ for every $\gamma\in SC(\alpha)$. \citep{Gebhardt2010};

\For{every $\gamma\in SC(\alpha)$} {

 Compute $d_{\gamma}$, the maximal dimension of a standard reduction simplex of $\gamma$ (\autoref{alg_standard_reduction_simplices});

}

Let $\beta$ be such that $d_{\beta}= \max\{d_\gamma;\ \gamma\in SC(\alpha)\}$;

Let $M_0$ be the intersection of all maximal standard reduction simplices of $\beta$;

Compute the period $p$ of the sequence $\{Z(\gamma^r)\}_{r>0}$ (Theorem~\ref{th_centralizers_are_periodic}).

\For{$i\in \{1,\ldots,p\}$}{

  Set $M_i=\varnothing$.

\If{$M_{i-1}\neq \varnothing$}{

  Compute generators $\{z_1,\ldots,z_m\}$ of the centralizer $Z(\gamma^i)$ (\cite{Franco2003b}).
 
  \While{$M_{i-1}\neq \varnothing$}{
 
   Extract $P$ from $M_{i-1}$, set $\mathcal S_1=\{P\}$, set $\mathcal S_2=\varnothing$. 
 
   \While{$\mathcal S_1\neq \varnothing$}{

    Extract $P$ from $\mathcal S_1$.
    
    \If{ $P^{z_1},\ldots,P^{z_m}$ all belong to $M_{i-1}$}{
      
      Set $\mathcal S_2=\mathcal S_2\cup\{P\}$.
      
      Add to $\mathcal S_1$ all elements in $\{P^{z_1},\ldots,P^{z_m}\}$ which are not already in $\mathcal S_2$.
      
    }
    \Else{
    
      Set $\mathcal S_1=\varnothing$ and $S_2=\varnothing$.
    }
   }
   \If{$\mathcal S_2\neq \varnothing$}{
      Set $M_i=M_i\cup \mathcal S_2$.
   }
  }
 }
}

\Return{$(M_p,{c_{\beta}})$;}

\end{algorithm}

We can explain \autoref{alg_computing_CRS_Artin} as follows. First we identify an element $\gamma\in SC(\alpha)$ whose canonical reduction system is standard, and we compute a standard simplex $M_0$ that contains $CRS(\gamma)$. Then we compute the period $p$ of the sequence $\{Z(\gamma^{r})\}_{r>0}$ and the sets $M_1,\ldots, M_{p}$. Each $M_i$ consists of the elements of~$M_0$ whose orbit under~$Z(\gamma^j)$ is contained in~$M_0$ for $j=1,\ldots,i$. Hence, once~$M_{i-1}$ is computed, the elements of~$M_i$ will be those $P\in M_{i-1}$ whose orbit under~$Z(\gamma^i)$ is contained in~$M_{i-1}$. We start to compute an orbit taking some $P\in M_{i-1}$, and iteratively computing the conjugate of each obtained element by $z_1,\ldots,z_m$ (the generators of~$Z(\gamma^i)$). The set~$\mathcal S_1$ consists of the elements still to be processed, while~$\mathcal S_2$ contains those that have already been processed without issue (i.e., their conjugates under $z_1,\ldots,z_m$ lie in~$M_{i-1}$). If some element outside~$M_{i-1}$ is found, all elements in $\mathcal S_2$ are discarded. If all computed elements belong to~$M_{i-1}$ and there are no more elements to process in~$\mathcal S_1$, then~$\mathcal S_2$ will contain an orbit under~$Z(\gamma^i)$ contained in~$M_{i-1}$, so all elements in~$\mathcal S_2$ are inserted into $M_i$. The algorithm continues as long as there are unprocessed elements in $M_{i-1}$. 

At the end of the process, for $i=1,\ldots,p$, the set $M_p$ will contain the elements of $M_0$ whose orbits under $\{Z(\gamma),\ldots,Z(\gamma^p)\}$ are contained in $M_0$ (so all elements in these orbits are standard). Since $\bigcup_{i=1}^{p}{Z(\gamma^i)}= \bigcup_{i>0}{Z(\gamma^i)}$, it follows from \autoref{th_standard_CRS_and_centralizers} that $M_p$ contains the set of essential irreducible parabolic subgroups for $\gamma$, hence $CRS(\gamma)=M_p$. The algorithm then returns $(M_p,c_\gamma)$, where $c_\gamma$ is a conjugating element from $\gamma$ to $\alpha$.

We will not undertake a study of the complexity of \autoref{alg_computing_CRS_Artin}, since it would depend on the cost of computing $p$ centralizers, and $p$ is a number that we do not know how to bound. In \autoref{sec_braids} we will provide an algorithm in the particular case of braid groups, and we will study its complexity.

\subsection{Some conjectures}

We saw in the previous subsection that finding a uniform upper bound for the period of $\{Z(\gamma^r)\}_{r>0}$ would improve the algorithm to compute canonical reduction systems (\autoref{conj_bound_for_period}). We will give now some conjectures which would allow us to find an alternative algorithm, without needing to compute any of the above centralizers. The conjectures we will state are true in the case of braid groups. 

\begin{conjecture}\label{conj_finite_orbits}
Let $G$ be an Artin-Tits group of spherical type, whose rank is $n$. For every $\alpha\in G$, the size of a finite orbit under $\alpha$ of a vertex in $\mathcal C(G)$ is bounded above by a number that depends only on $n$.
\end{conjecture}

\begin{proposition}\label{prop_finite_orbits_in_braids}
 \autoref{conj_finite_orbits} holds if $G$ is a braid group. Every finite orbit of an element $\alpha\in \mathcal B_{n+1}$ (or rank $n$) has size at most $n+1$.
\end{proposition}
\begin{proof}
In the case of a braid group of rank $n$ (with $n+1$ strands), a finite orbit of an element $\alpha$ can be of two types. If it consists of essential curves for $\alpha$, it is a simplex of $\mathcal C(G)$ whose vertices do not admit any inclusion relation. Hence, it consists of at most $(n+1)/2$ vertices. 

If a finite orbit does not consist of essential curves, each curve belongs to a periodic component of the decomposition of $\alpha$ along its canonical reduction system (and each component encloses at least 3 punctures, otherwise it could not contain non-essential curves). Let $r<(n+1)/3$ be the number of curves in $CRS(\alpha)$ which enclose the curves of the studied orbit (these curves of $CRS(\alpha)$ form themselves an orbit of $\alpha$). Then $\alpha^r$ preserves each of these curves from $CRS(\alpha)$, and the restriction of $\alpha^r$ to the subsurface corresponding to each of these curves is periodic, with $m\leq (n+1)/r$ punctures. Hence, either $\alpha^{r(m-1)}$ or $\alpha^{rm}$ preserves each of these curves, and the restriction is a pure periodic braid. Since a pure periodic braid preserves all curves, either $\alpha^{r(m-1)}$ or $\alpha^{rm}$ fixes all curves of the original finite orbit. Hence, the original finite orbit has size at most $rm\leq n+1$.
\end{proof}

The bound given by the above result is sharp, as we can consider the example of the braid $\alpha=\sigma_n\sigma_{n-1}\cdots \sigma_1\in \mathcal B_{n+1}$ and the standard curve $C$ enclosing the first two punctures. The curve $C^{\alpha^i}$ will enclose the punctures $i+1$ and $i+2$ (module $n+1$), so $C\neq C^{\alpha^i}$ for $i=1,\ldots,n$. But $\alpha^{n+1}=\Delta^2$, which is central, so $C^{\alpha^{n+1}}=C$. Therefore, the orbit of $C$ under $\alpha$ has size $n+1$, which is the bound given in \autoref{prop_finite_orbits_in_braids}.

This leads us to the following:

\begin{conjecture}\label{conj_bound_on_finite_orbits}
The bound conjectured in \autoref{conj_finite_orbits} is the maximal order of a periodic element in $G$.
\end{conjecture}

If \autoref{conj_bound_on_finite_orbits} is true (or if \autoref{conj_finite_orbits} is true and we know an explicit bound), then a vertex of $\mathcal C(G)$ with a finite orbit under $\alpha$ can be detected as follows: it would be a fixed vertex for $\alpha^r$ for some bounded $r$, and it would be seen as a \emph{standard} fixed vertex for some element in~$SC(\alpha^r)$. Since $r$ is bounded and the number of standard vertices is bounded too, we have only a finite number of vertices to check.

Suppose that $M_0$ is the standard reduction simplex of $\gamma\in SC(\alpha)$, containing $CRS(\gamma)$. Let~$\delta$ be a conjugate of~$\alpha$ admitting a standard vertex~$Q$ of whose orbit under~$\delta$ is finite. And let $c$ be a conjugating element such that $\gamma^c=\delta$. Then, for every $P\in M_0$, we can check whether~$P^c$ is adjacent to $Q$. If~$P$ belongs to $CRS(\gamma)$, then~$P^c$ belongs to $CRS(\delta)$, hence~$P^c$ and~$Q$ are adjacent. Therefore, if~$P^c$ and~$Q$ are not adjacent, we can discard $P$ (and all its images under elements of $Z(\gamma)$) as it does not belong to $CRS(\gamma)$.

Since all elements with finite orbit under $\alpha$ can be obtained as a subgroup $Q$ satisfying the condition of the previous paragraph, it follows that the elements $P\in M_0$ which are not discarded will belong to $CRS(\gamma)$. Hence, we would have an algorithm to compute $CRS(\gamma)$, and then $CRS(\alpha)$.

\section[The case of braid groups: detecting curves from the CRS]{The case of braid groups: detecting curves from the canonical reduction system }\label{sec_braids}

In this section we will consider $G=\mathcal B_n$, the braid group on $n$ strands,that is, the Artin-Tits group $A_{n-1}$. All previously shown results are valid in this case, applied to the curve complex~$\mathcal C(\discn)$. So we will consider curves instead of irreducible parabolic subgroups.

We will give an algorithm to compute $CRS(\alpha)$ for a braid $\alpha\in \mathcal B_n$. We already know how to compute a standard reduction simplex $M_0$ for a braid $\gamma\in SC(\alpha)$, which contains $CRS(\gamma)$. We need to determine, for each curve in $M_0$, whether it belongs to $CRS(\gamma)$ or not. For simplicity, we will assume that $\gamma=\alpha$, so $CRS(\alpha)$ is standard and is contained in a standard multicurve.

\subsection{Canonical reduction system and centralizers}

Let us see how we can improve \autoref{th_standard_CRS_and_centralizers} in the case of $\mathcal B_n$, relating the centralizer of a braid~$\alpha$ with the curves in $CRS(\alpha)$. This is not the fastest way to detect curves from $CRS(\alpha)$, but it seems to be a promising method which could eventually be generalized to other Artin-Tits groups of spherical type.

We distinguish between the case where $\alpha$ is pure ---since it allows for a more efficient approach--- and the general case.

\begin{definition}\label{def_subsurface}
Let $\alpha$ be a braid with standard canonical reduction system. Let $M$ be a standard multicurve preserved by $\alpha$. Notice that $M_+=M\cup CRS(\alpha)$ is also a standard multicurve preserved by $\alpha$, and that each connected component of $\discn\setminus (M_+ \cup \partial(\discn))$ is an open punctured disk. If $C$ is either $\partial(\discn)$ or a curve of $M$, we define~$D_C$ to be $X_C\cup C$, where~$X_C$ is the connected component of $\discn\setminus (M_+ \cup \partial(\discn))$ which is the outermost component enclosed by~$C$. 
\end{definition}

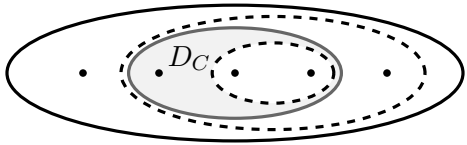
\begin{figure}[H]
    \centering
    \begin{tikzpicture}

\draw[very thick] (2,0) ellipse (3 and 0.9);
\fill (0,0) circle (0.05);
\fill (1,0) circle (0.05);
\fill (2,0) circle (0.05);
\fill (3,0) circle (0.05);
\fill (4,0) circle (0.05);

\draw[dashed, very thick] (2.5,0) ellipse (2 and 0.75);

\filldraw[color=black!60, fill=black!5, very thick] (2,0) ellipse (1.4 and 0.6);
\fill[black] (1,0) circle (0.05);
\fill[black] (2,0) circle (0.05);
\fill[black] (3,0) circle (0.05);

\filldraw[dashed, fill=white, very thick] (2.5,0) ellipse (0.8 and 0.4);
\fill (2,0) circle (0.05);
\fill (3,0) circle (0.05);

\draw (1.4,0.2) node{} node{$D_C$};
\end{tikzpicture}
    \caption{The $5$-punctured disk with an example of a CRS (dashed) and  of $D_C$ (the gray area and the gray curve $C$). In this case $D_C$ is homeomorphic to the $2$-punctured disk simply by collapsing to a puncture the region enclosed in the curve in the innermost dashed curve.}
    \label{fig:DC}
\end{figure}

Notice that $D_C$ is homeomorphic to a $r$-punctured disk, whose boundary component is $C$ (see \autoref{fig:DC}).

\begin{definition}\label{def_restriction}
Let $\alpha$ be a braid, let $M$ be a standard multicurve preserved by $\alpha$, and let $C$ be either $\partial(\discn)$ or a curve of $M$. Let $m$ be the smallest positive integer such that~$\alpha^m$ fixes~$C$. Then we define the {\em component of~$\alpha$ corresponding to $C$}, denoted $\alpha_C$, to be the restriction of~$\alpha^m$ to~$D_C$. 
\end{definition}

\begin{figure}[h]\centering
\begin{tikzpicture}

    \draw[line width=4pt,white] (2.5,-2) ellipse (1cm and .4cm); 
    \draw[line width=4pt,white] (2.5,-2) ellipse (2.2cm and .65cm); 
    \draw[dashed,thick] (2.5,-2) ellipse (2.2cm and .65cm); 
    \filldraw[black!60,fill=black!5,thick] (2,-2) ellipse (1.6cm and .5cm); 
    \filldraw[dashed,fill=white,thick] (2.5,-2) ellipse (1cm and .4cm); 

    \draw[] (0.5,0) to[out=-90,in=90] (1.5,-2);
    \thirdpartialcurve{0.75}{0.77}
    \draw[] (2.5,0) to[out=-90,in=90] (3.5,-2);
    \firstpartialcurve{0.77};
    \draw[black] (3,0) to[out=-90,in=90] (1,-2);
 
    \filldraw[dashed,fill=white] (1.5,0) ellipse (1cm and .4cm); 

    \draw[black] (1,2) to[out=-90,in=90] (3,0);
    \secondpartialcurve{0.27}{0.74}
    \draw[line width=4pt,white] (4,2) to [out=-90,in=90] (4,0);
    
    \draw[] (3.5,2) to[out=-90,in=90] (2.5,0);
    \draw[] (1.5,2) to[out=-90,in=90] (0.5,0);
    \draw[line width=4pt,white] (2.5,2) ellipse (1cm and .4cm); 
    \filldraw[black!60,fill=black!5,thick] (2,2) ellipse (1.6cm and .5cm); 
    \filldraw[dashed,fill=white,thick] (2.5,2) ellipse (1cm and .4cm); 

   \draw[line width=2pt,white] (2.5,2) ellipse (2.2cm and .65cm); 
    \draw[dashed,thick] (2.5,2) ellipse (2.2cm and .65cm); 

    \draw[dashed,thick] (.3,2) to (.3,-2) (4.7,2) to (4.7,-2);

    \draw[black] (1,2) node[circle,inner sep=1pt,fill]{} node[left]{};
    \draw[] (2,2) node[circle,inner sep=1pt,fill]{} node[left]{};
    \draw[] (3,2) node[circle,inner sep=1pt,fill]{} node[left]{};
    \draw[black] (4,2) node[circle,inner sep=1pt,fill]{} node[left]{};

    \draw[] (1,-2) node[circle,inner sep=1pt,fill]{} node[left]{};
    \draw[black] (2,-2) node[circle,inner sep=1pt,fill]{} node[left]{};
    \draw[black] (3,-2) node[circle,inner sep=1pt,fill]{} node[left]{};
    \draw[] (4,-2) node[circle,inner sep=1pt,fill]{} node[left]{};
    
\end{tikzpicture}
\caption{A possible component of $\alpha$ with respect to a curve $C$ (solid circle). In this case $\alpha_C$ is the full twist on $2$-strands.}
\label{fig:DC-component}
\end{figure}
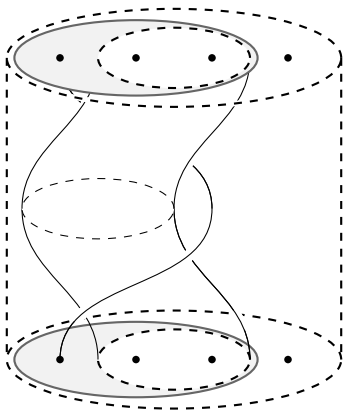

Notice that, if $D_C$ is homeomorphic to a $r$-punctured disk, the component $\alpha_C$ can be seen as a braid on $r$ strands (see \autoref{fig:DC-component}).

Let us improve \autoref{th_standard_CRS_and_centralizers} in the case of pure braids.

\begin{theorem}\label{pure} Let $\alpha\in \mathcal{B}_n$ be a pure braid with standard canonical reduction system. A standard curve~$C$ belongs to $CRS(\alpha)$ if and only $C^z$ is standard for every $z\in Z(\alpha)$.
\end{theorem}

\begin{proof}
First, for every element $z\in Z(\alpha)$, we have $CRS(\alpha)^z = CRS(\alpha^z)=CRS(\alpha)$, which is a set of standard curves. Hence, if $C\in CRS(\alpha)$, then~$C^z$ will be standard.

Now suppose that $C\notin CRS(\alpha)$. If~$C$ has non-trivial intersection with some curve $C'\in CRS(\alpha)$, then the orbit of $C$ under $\alpha$ is infinite. Since the set of standard curves is finite, it follows that $C^{\alpha^m}$ will be non-standard for some $m>0$. Since $\alpha^m\in Z(\alpha)$, we have the desired implication for this case.

We can then assume that $C$ has trivial intersection with all curves in $CRS(\alpha)$. This means that $C$ is contained in~$D_{C'}$ for some $C'\in CRS(\alpha)\cup \{\partial(\discn)\}$. 

If~$\alpha_{C'}$ is pseudo-Anosov, no power of~$\alpha$ fixes~$C$. Hence, as before, $C^{\alpha^m}$~is not standard for some $m>0$.

If~$\alpha_{C'}$ is periodic, since $\alpha$ is pure, then $\alpha_{C'}$ is a power of the full twist in $D_{C'}$. Hence, any braid with support on~$D_{C'}$ commutes with $\alpha$. We can now take as~$z\in \mathbb{Z}(\alpha)$ any braid with support on~$D_{C'}$ that sends~$C$ to a non-standard curve (this braid exists because~$D_{C'}$ has at least three punctures, since it can contain a non-degenerate curve~$C$). 
\end{proof}

\autoref{pure} is sufficient to obtain an algorithm which computes $CRS(\alpha)$ for any $\alpha$, since $CRS(\alpha)=CRS(\alpha^m)$ for every $m>0$, and it suffices to take $m$ such that $\alpha^m$ is pure. But the number $m$ could be exponentially big with respect to the number $n$ of strands, as shown in the following example: Let $p_1,p_2,\ldots,p_r$ be the first $r$ prime numbers. Let $\alpha\in \mathcal B_n$ with $n=p_1+p_2+\cdots+p_r$, inducing a permutation with a cycle of length $p_i$ for $i=1,\ldots,r$. Then, the smallest $m>0$ such that $\alpha^m$ is pure is $m=p_1p_2\cdots p_r$. As $r$ tends to infinity, $m$ grows exponentially with respect to $n$.

Hence, we will give more precise results that will provide a faster algorithm to compute $CRS(\alpha)$. We start with some basic results on permutations.

\begin{remark}\label{remark_permutation}
Let $\pi\in Sym_n$ be a permutation of $\{1,\ldots,n\}$. If $\pi$ is a power of a cycle, $\pi=(a_1 \cdots a_r)^l$, then $\pi$ is the identity if and only if $\pi(a_i)=a_i$ for some $i$.   
\end{remark}

Let $o_\alpha(C)=\{C^{\alpha^r}\,|\, r>0\}$ be the  (positive) orbit of a curve $C$ under the action of a braid~$\alpha$. Similarly, every braid induces a permutation on the punctures of $\discn$, so we can define $o_\alpha(x)=\{\alpha^{r}(x)\,|\, r>0\}$ to be the orbit of a puncture $x$ under the action of $\alpha$. This latter orbit is always finite, smaller than or equal to $n$.

Given a curve $C$ whose orbit under $\alpha$ is a multicurve, we can consider $D_C$, and we will also consider the orbits of the punctures of $D_C$. It is important to notice that a puncture $x$ of $D_C$ may correspond to either a puncture of $\discn$ or a curve in $CRS(\alpha)$. In any case, its orbit under~$\alpha$ will be finite, and we will denote it by $o_\alpha(x)$, as above.

\begin{lemma}\label{lemma_power}
Let $\alpha\in \mathcal{B}_n$ be a braid with a standard canonical reduction system. Let $M$ be a standard multicurve such that $M=o_{\alpha}(C)$ for some standard curve $C$. Suppose that $\alpha_C$ is periodic, and let $x$ and $y$ be two punctures of $D_C$. Then $l:=\lcm(|o_\alpha(x)|,|o_\alpha(y)|)\leq n$, the braid~$\alpha^l$ fixes~$C$ and $(\alpha^l)_C$ is pure (so it is a power of the full twist).
\end{lemma}
 
\begin{proof}
 Let $m_1$ be the number of curves in $M$. Since $M$ is an orbit, the braid $\alpha^{m_1}$ fixes $C$, and the restriction of $\alpha^{m_1}$ to $D_C$ is precisely $\alpha_C$, which is periodic. Moreover, if $r$ is the number of punctures of $D_C$, we have $r\leq n/m_1$. Then, either $(\alpha^{m_1})^{r-1}$ or $(\alpha^{m_1})^r$ preserves $C$ and its restriction to $D_C$ is pure due to Brouwer-Kérékjartó-Eilenberg Theorem \citep{Constantin1994}. Hence, some power $\alpha^{l'}$ with $1\leq l'\leq n$ fixes $C$ and all punctures of $D_C$.

 Now, let $x$ and $y$ be two punctures of $D_C$, and let $l=\lcm(|o_\alpha(x)|,|o_\alpha(y)|)$. The number $l$ is the smallest exponent such that $\alpha^l$ fixes $x$ and $y$, therefore $l\leq l'\leq n$. Now $x$ is enclosed by $C$, hence $\alpha^l(x)=x$ is enclosed by $C^{\alpha^l}$. But the orbit of the curve $C$ is a multicurve (and no two curves can be nested), hence $C^{\alpha^l}=C$. It follows that $\alpha^l$ fixes $C$ and two punctures~$x$ and~$y$ of~$D_C$.  Since the permutation corresponding to a periodic braid on $r$ strands is the power of either a cycle of length $r$ or a cycle of length $r-1$, it follows from \autoref{remark_permutation} that the permutation induced by $\alpha^l$ on the punctures of $D_C$ is trivial. That is, $(\alpha^l)_C$ is pure. 
\end{proof}

We can now see that it is not necessary to raise $\alpha$ to an exponentially big power, to check whether a given standard curve belongs to $CRS(\alpha)$.

\begin{theorem}\label{non-pure} Let $\alpha\in \mathcal{B}_n$ be a braid with a standard canonical reduction system. Let $M$ be a standard multicurve such that $M=o_{\alpha}(C)$ for some standard curve $C$. Let $x$ and $y$ be two punctures of $D_C$, and let $l=\lcm(|o_\alpha(x)|,|o_\alpha(y)|)$. Then $C$ belongs to $CRS(\alpha)$ (hence $M$ is contained in $CRS(\alpha)$) if and only if either $l>n$ or $C^z$ is standard for every $z\in Z(\alpha^l)$.
\end{theorem}

\begin{proof}
Suppose that $C\in CRS(\alpha)$. Then $C\in CRS(\alpha^l)$ for every $l>0$. Now, if $z\in Z(\alpha^l)$ then  $C^z\in CRS(\alpha^l)^{z}=CRS((\alpha^l)^z)=CRS(\alpha^l)=CRS(\alpha)$. Since this is a set of standard curves, it follows that $C^z$ is standard.

Now assume that~$C$ does not belong to $CRS(\alpha)$. Then~$C$ must belong to a periodic component of~$\alpha$. That is, there is some curve $C'\in CRS(\alpha)\cup \{\partial(\discn)\}$ such that~$C$ is a non-degenerate curve in $D_{C'}$. Notice that, in this case, $\alpha_C$~is periodic. Moreover, since the orbit of $C$ is finite, no curve of $CRS(\alpha)$ can intersect~$C$, so the punctures of~$D_C$ are punctures of~$D_{C'}$. 

Let $x$, $y$ be two punctures of $D_{C}$ (which are also punctures of $D_{C'}$). Applying \autoref{lemma_power} to the orbit of $C'$ under $\alpha$, we know that  $l=\lcm(|o_\alpha(x)|,|o_\alpha(y)|)\leq n$ and that $(\alpha^l)_{C'}$ is a power of a full twist on at least three strands (because the punctured disk $D_{C'}$ admits $C$ as a non-degenerate curve). Now we can consider any braid $\beta$ with support in $D_{C'}$ which sends $C$ to a non-standard curve. Then $\beta$ commutes with $\alpha^l$ (since $\alpha^l$ fixes $C'$ and the component $(\alpha^l)_{C'}$ is a full twist). Therefore, if $C\notin CRS(\alpha)$ then $l\leq n$ and there exists $\beta\in Z(\alpha^l)$ such that $C^\beta$ is not standard.
\end{proof}

We can use the latter result to develop an algorithm that computes $CRS(\alpha)$ for any braid $\alpha\in \mathcal B_n$. But in the next subsection we will see a much faster approach. In any case, one could possibly extend these results (and not the ones in the next subsection) to other Artin-Tits groups of spherical type. And we think that the results in this section, relating elements of $CRS(\alpha)$ to centralizers of powers of $\alpha$, have their own interest.

\subsection{Gluing components}\label{ss_gluing_components}

Suppose that $\alpha\in \mathcal B_n$ and that $M$ is a standard reduction simplex for $\alpha$ that contains $CRS(\alpha)$. We will see a fast algorithm to determine which curves of $M$ belong to $CRS(\alpha)$.

Given a curve $C\in M$, in \autoref{def_subsurface} we defined the component $X_C$. Since $M$ contains $CRS(\alpha)$, we have that $M^+=M$, and that $X_C$ is just a component of $\discn\setminus(M\cup\{\partial(\discn)\}$, namely, the one having $C$ as its external boundary. Then $D_C=X_C\cup C$, and the braid $\alpha_C$ is the restriction of $\alpha^m$ to $D_C$, where $m$ is the smallest positive integer such that $\alpha^m$ fixes $C$ (\autoref{def_restriction}). Since we will apply these definitions with distinct sets of curves, we will use, instead of $D_C$ and $\alpha_C$, the notation $D_{M,C}$ and $\alpha_{M,C}$ to indicate that $M$ is the multicurve used to decompose $\discn$.

Given $\alpha$, computing $\alpha_{M,C}$ is very fast, following these steps:
\begin{enumerate}

 \item Consider the punctures $\{x_1,\ldots,x_r\}$ of $D_C$. If $x_i$ is a puncture of $\discn$, let $y_i=x_i$. If $x_i$ is a curve, choose a puncture $y_i$ of $\discn$ enclosed by $x_i$. Then $\{y_1,\ldots,y_r\}$ is a set of punctures of $\discn$.

 \item Compute $m$, considering the irreducible parabolic subgroup $P$ corresponding to $C$, denoting $z_{(0)}=z_P$, and iteratively computing $z_{(i)}=(z_{(i-1)})^{\alpha}$ until the first repetition $z_{(m)}=z_{(0)}$. Each $z_{(i)}$ will be a positive element, as $M$ is standard.

 \item Compute $\alpha^m$, and keep only the strands starting at $\{y_1,\ldots,y_r\}$. That is, from any word representing $\alpha^m$, keep only the letters (strands' crossings) corresponding to two strands of $\{y_1,\ldots,y_r\}$, with the appropriate index shifting. 
 
\end{enumerate}

The resulting braid will be $\alpha_{M,C}$. The process of deleting strands can be done to each simple element (permutation braid), considering the permutation induced on the strands we want to keep. This has complexity $O(n)$ for each simple element. Since $m\leq n/2$, and the complexity of computing a normal form of a braid which is the product of $l$ simple elements is $O(l^2n\log n)$, the complexity of the whole procedure is $O((Nn)^2n\log n)=O(N^2n^3\log n)$, where $N$ is the canonical length of $\alpha$. This is the same complexity of computing the normal form of $\alpha^m$.

Now, the curve $C\in M$ is a boundary of exactly two subsurfaces: it is the external boundary of $D_{M,C}$ and it is a puncture of $D_{M,C'}$, where $C'\in M\cup\{\partial(\discn)\}$ is the smallest curve in this set that encloses $C$. 

\begin{proposition}
With the above conditions, $C\notin CRS(\alpha)$ if and only if $\alpha_{M\setminus \{C\},C'}$ is periodic.
\end{proposition}

\begin{proof}
Since $M$ contains $CRS(\alpha)$, $M$ is an adequate reduction system of $\alpha$, that is, each component of $\alpha$ determined by $M$ is either periodic or pseudo-Anosov. 

If $\alpha_{M\setminus \{C\},C'}$ is periodic, this means that $M\setminus\{C\}$ is also adequate, as $D_{M\setminus \{C\},C'}$, obtained by gluing $D_{M,C}$ and $D_{M,C'}$ is the only new subsurface of $\discn$ that we obtain when we remove $C$ from $M$.

Since $CRS(\alpha)$ is the minimal adequate reduction system (by inclusion), this implies that $C\notin CRS(\alpha)$. 

Conversely, suppose that $C\notin CRS(\alpha)$. Since $C$ is preserved by a power of $\alpha$, it cannot be contained into a pseudo-Anosov component of $\alpha$ with respect to $CRS(\alpha)$. Hence, it is contained into a periodic component, whose associated subsurface is $D_{CRS(\alpha),C''}$. That subsurface is subdivided into smaller subsurfaces by the curves of $M$ which it contains. But the component of $\alpha$ associated to each of them must be periodic (some power will become a full twist).  This happens also if we remove the curve $C$ from $M$. Hence, the component $\alpha_{M\setminus \{C\},C'}$ is periodic.
\end{proof}

As a conclusion, we have a very fast algorithm to determine whether a curve of $M$ belongs to $CRS(\alpha)$ or not. We just need to compute $\beta=\alpha_{M\setminus \{C\},C'}$ (say it has $r$ strands), and then check whether it is periodic, which will happen if and only if either $\beta^{r-1}$ or $\beta^r$ is a full twist on $r$ strands, that is, if its left normal form is a power of $\Delta$ (on $r$ strands). 

Notice that the power to which we had to raise $\alpha$ multiplied the power to which we need to raise $\beta$ gives a number that is at most $n$. Hence, the complexity of the whole procedure is the same as computing the normal form of $\alpha^m$, for some $m\leq n$, hence it is $O(N^2n^3\log n)$.

\section[An algorithm to compute the CRS of a braid]{An algorithm to compute the canonical reduction system of a braid}

We can finally provide an algorithm to compute $CRS(\alpha)$ for any braid $\alpha$, using the procedures explained in the previous sections.

\bigskip 

\begin{algorithm}[H]\caption{For computing the canonical reduction system of a braid.}\label{alg_computing_CRS_braids}

\SetKwInOut{Input}{Input}\SetKwInOut{Output}{Output}

\BlankLine

\Input{An integer $n\geq 2$ and a braid $\alpha\in\mathcal{B}_n$.}
\Output{$CRS(\alpha)=\mathcal{C}^c$ as a pair $(\mathcal C,c)$, where $\mathcal{C}$ is a standard simplex and $c\in \mathcal B_n$.}
\BlankLine
\BlankLine

Compute $SC(\alpha)$, and a conjugating element $c_\gamma$ such that $\gamma^{c_{\gamma}}=\alpha$ for every $\gamma\in SC(\alpha)$. \citep{Gebhardt2010};

\For{every $\gamma\in SC(\alpha)$} {

 Compute $d_{\gamma}$, the maximal dimension of a standard reduction simplex of $\gamma$ (\autoref{alg_standard_reduction_simplices});

}

Let $\beta$ be such that $d_{\beta}= \max\{d_\gamma;\ \gamma\in SC(\alpha)\}$;

Let $M$ be the intersection of all maximal standard reduction simplices of $\beta$;

\For{every $C\in M$}{

 Compute $\delta=\beta_{M\setminus\{C\},C'}$ (\autoref{ss_gluing_components});
 
 \If{$\delta$ is periodic}{
   Set $M=M\setminus o_\beta(C)$;
 }
}
\Return{$(M,{c_{\beta}})$;}

\end{algorithm}

\bigskip

Let us study the complexity of \autoref{alg_computing_CRS_braids}.

The complexity of computing $CRS(\alpha)$ is generically polynomial with respect to $n$ and to the length of $\alpha$, but this is due to the fact that a generic braid is pseudo-Anosov (thus its canonical reduction system is trivial). The worst case complexity of computing $CRS(\alpha)$ depends on the following data. We assume that $\alpha$ is given as a product of $\ell$ simple elements. 

\begin{itemize}

\item The number $K_\alpha := |SC(\alpha)|$. Again, generically~$K_\alpha$ is very small, but its maximal possible value is unknown in general. If $N$ is the canonical length of $\alpha$, an upper bound (very far from optimal) is $(n!-2)\ell$, as every braid in $SC(\alpha)$ can be written as $\Delta^k \beta_1\cdots \beta_{N'}$ with $k$ and $N'$ fixed, and $N'\leq N$. Since the number of possible proper simple factors is $n!-2$, we obtain the bound.

\item The number $T_\alpha$ of cyclic sliding operations starting from $\alpha$ until reaching the first repetition, thus completing a sliding circuit. Although this number is usually very similar to the length of $\alpha$, its maximal value is not known. A (very far from optimal) bound is given by ${n\choose 2}\ell + (n!-2)\ell$, since one needs at most ${n\choose 2}\ell$ cyclic slidings to send $\alpha$ to its Super Summit Set, and all the obtained elements from there on have the same canonical length.

\item The cost of computing the preferred prefix $\mathfrak{p}(\alpha)$, which is $O(n\log(n))$ \citep{Gebhardt2010}.

\item The cost of computing the left normal form of a braid of canonical length $\ell$, which is $O(\ell^2 n\log(n))$ \citep{Epstein1992}.

\end{itemize}

Starting with the braid $\alpha$, we perform $T_\alpha$ cyclic slidings. If an element $\beta$ is already in normal form (as a product of $\ell$ simple factors), recomputing the normal form of $\mathfrak{p}(\beta)^{-1}\beta\mathfrak{p}(\beta)$ costs $O(\ell n\log(n))$. Therefore, applying $T_\alpha$ cyclic slidings to $\alpha$ can be done in time $O(T_\alpha \ell n\log(n))$.

According to \citep{Gebhardt2010}, to construct the whole set $SC(\beta)$, for every $\gamma\in SC(\beta)$ we must compute all minimal simple elements conjugating $\gamma$ to an element in $SC(\beta)$. We then add to~$SC(\beta)$ the new elements (the conjugates of~$\gamma$ by these minimal elements) and continue the process until completion.  Let~$L_\beta$ be the cost of computing such minimal elements for a given element in $SC(\beta)$, that again is generically small. Then, the complexity of computing the whole $SC(\alpha)$ is $O(T_\alpha \ell n\log(n) + L_\alpha K_\alpha)$.

\autoref{alg_computing_CRS_braids} continues by computing maximal standard reduction simplices for each $\gamma\in SC(\alpha)$, in order to obtain the maximal possible $d_\gamma$. The complexity of computing such simplices was studied after \autoref{alg_standard_reduction_simplices}, and it is $O(\ell n^6+3^{n/3})$. Hence, the computation of $\beta$ and $d_{\beta}$ is done in $O(K_\alpha (\ell n^6+3^{n/3}))$.

Once $\beta$ is found, we say in the previous section that determining which curves of $M$ belong to $CRS(\beta)$ has complexity $O(\ell^2n^3\log n)$.

Joining all the above ingredients, we have shown the following:

\begin{theorem}
If $\alpha\in \mathcal B_n$ has canonical length $\ell$, the complexity of computing $CRS(\alpha)$ using \autoref{alg_computing_CRS_braids} is
$$
O\left(T_\alpha \ell n\log(n) +
K_\alpha (L_\alpha+\ell n^6+3^{n/3})+\ell^2n^3\log n\right).
$$
\end{theorem}

As we said above, the number $T_\alpha$, $L_\alpha$ and $K_\alpha$ are generically very small, and one could a priori take advantage of the properties of canonical reduction systems to improve the number~$3^{n/3}$.

We prefer to leave \autoref{alg_computing_CRS_braids} as it is, for simplicity.
\bigskip

\noindent{\textbf{\Large{Acknowledgments}}}

We are deeply grateful to Ivan Marin, for his help in proving \autoref{th_centralizers_are_periodic}, and for very useful discussions.

The first and third authors were supported by the research project PID2022-138719NA-I00, and the second author by the project PID2024-157173NB-I00. The first author was also financed by an individual Ram\'on y Cajal 2021 Fellowship RYC2021-032540-I. All these projects were funded by the Spanish Ministry of Science and Innovation\\ (MCIN/AEI/10.13039/501100011033/FEDER, UE).

The third author also acknowledges support from the Grant QUALIFICA by Junta de Andalucía grant number QUAL21 005 USE, from the Swiss Government Excellence Scholarship and from Swiss NSF grant 200020-200400. 
The third author is a member of the Gruppo Nazionale per le Strutture Algebriche, Geometriche e le loro Applicazioni (GNSAGA) of the Istituto Nazionale di Alta Matematica (INdAM). He also thanks the Departamento de Álgebra of the Universidad de Sevilla for the kind hospitality.

\medskip

\medskip
\bibliography{main}

\bigskip\bigskip{\footnotesize%

\noindent
\textit{\textbf{María Cumplido} \\ 
Instituto de Matemáticas de la Universidad de Sevilla (IMUS) and Departamento de Álgebra, Universidad de Sevilla, Spain.} \par
 \textit{E-mail address:} \texttt{\href{mailto:cumplido@us.es}{cumplido@us.es}}

\noindent
\textit{\textbf{Juan González-Meneses} \\ 
Instituto de Matemáticas de la Universidad de Sevilla (IMUS) and Departamento de Álgebra, Universidad de Sevilla, Spain.} \par
 \textit{E-mail address:} \texttt{\href{mailto:meneses@us.es}{meneses@us.es}}

 \noindent
\textit{\textbf{Davide Perego} \\ 
Section de Mathématiques, Université de Genève, Switzerland.} \par
 \textit{E-mail address:} \texttt{\href{mailto:dperego9@gmail.com}{dperego9@gmail.com}}

\end{document}